\documentclass[10pt]{amsart}  
\usepackage{epsfig}
\usepackage{color}
\usepackage{amsmath, amssymb, euscript,  enumerate}
\usepackage{verbatim}
 \usepackage{pxfonts}
\usepackage{graphicx}
\usepackage{subfigure}
\usepackage{color} 
\newcommand{\nc}{\newcommand}

\newtheorem{theorem}{Theorem}[section]

\newtheorem{corollary}[theorem]{Corollary}

\newtheorem{lemma}[theorem]{Lemma}

\newtheorem{proposition}[theorem]{Proposition}
\newtheorem{remark}[theorem]{Remark}

\newcommand{\beq}{\begin{equation}}
\newcommand{\eeq}{\end{equation}}
\newcommand{\bcl}{\begin{center}}
\newcommand{\ecl}{\end{center}}
\newcommand{\re}{\mathbb{R}}

\newcommand{\cC}{{\mathcal C_b}}
\newcommand{\cD}{{\mathcal D}}
\newcommand{\cS}{{\mathcal S}}

\newcommand{\cM}{{\mathcal M}}
\newcommand{\cA}{{\mathcal A}}
\newcommand{\cP}{{\mathcal P}}

\newcommand{\e}{\varepsilon}

\newcommand{\p}{\partial}
\newcommand{\La}{\triangle}
\newcommand{\D}{\nabla}
\newcommand{\ra}{\rightarrow}
\newcommand{\pa}{\partial}
\newcommand{\supp}{\operatorname{supp}}
\newcommand{\R}{{\mathbb R}}
\newcommand{\mS}{{\mathbb S}}
\newcommand{\vp}{\varphi}
\newcommand{\oa}{\overline{\alpha}}

\newcommand{\al}{\alpha}
\newcommand{\be}{\beta}
\nc \sms {\smallskip}
\nc{\fp}{\noindent}
\nc{\qd}{\qquad\qquad}

\title[Asymptotic Behavior and Nonlinear Eigenvalue Problems]
{Asymptotic Behavior in Degenerate Parabolic Fully Nonlinear equations
 and its application to Elliptic Eigenvalue Problems
}

\author[Soojung Kim]{Soojung Kim}
\address{Soojung Kim :
School of Mathematical Sciences,
Seoul National University,
1 Gwanrk-ro, Gwanak-Gu, Seoul 151-747, South Korea}
\email{soojung26@gmail.com }

\author[Ki-ahm Lee]{Ki-ahm Lee}
\address{Ki-ahm Lee:
School of Mathematical Sciences,
Seoul National University,
1 Gwanrk-ro, Gwanak-Gu, Seoul 151-747, South Korea\& Korea Institute for Advanced Study, Seoul,130-722, Korea}
\email{kiahm@math.snu.ac.kr}

\begin{document}

\begin{abstract}
\noindent  We study the asymptotic behavior of the nonlinear
parabolic flows $u_{t}=F(D^2 u^m)$  when $t\ra \infty$ for $m\geq 1$, and   the geometric properties
for solutions of the following elliptic nonlinear eigenvalue problems:
\begin{equation*}
\begin{split}
F(D^2 \vp) &+ \mu\vp^{p}=0, \quad \vp>0\quad\text{in $\Omega$}\\
\vp&=0\quad\text{on $\p\Omega$}
\end{split}
\end{equation*}
posed in a (strictly) convex and smooth domain $\Omega\subset\re^n$ for $0< p \leq 1,$  where $F(\cdot)$ is uniformly elliptic, positively homogeneous of order one and concave. 
We establish that  $\log (\vp)$ is  concave in the case $p=1$ and that the function $\vp^{\frac{1-p}{2}}$ is
 concave for  $0<p<1.$\end{abstract}
\maketitle

\section{Introduction}\label{sec-intro}
In this paper, we consider the asymptotic behavior of $u$ satisfying
\begin{equation}\label{eq-intro-1}
\left\{
\begin{array}{ll}
u_t(x,t)-F(D^2 u^m(x,t))=0\qquad &\text{in $Q=\Omega\times (0,+\infty)$},\\
u(x,0)>0\quad&\text{in $\Omega,$}\\
u(x,t)=0\quad &\text{on $\in \partial\Omega\times(0,+\infty),$}
\end{array}\right. \qquad \qquad
\end{equation}
and then we show a renormalized flow converges to $\vp(x)$ which satisfies the following nonlinear eigenvalue problem
\begin{equation}\label{eq-intro}
\left\{
\begin{array}{ll}
 F(D^2 \vp)+  \mu\vp^{p} &=0,\quad\text{in $\Omega,$}\\
  \vp&>0\quad\text{in $\Omega,$}\\
 \vp&=0\quad\text{on $\p\Omega$}
\end{array}\right. \qquad \qquad
\end{equation}
for some $\mu>0.$

Such parabolic approach to nonlinear eigenvalue problem has been considered at \cite{LV2} for Laplace operator and extended to Fully nonlinear operator at \cite{KsL} with super-linear exponent.(i.e. $1<p<\textsf{p}_{\Omega,F}$ for some critical number $\textsf{p}_{\Omega,F}>1$.)
 In this paper, we consider linear and sublinear case ($0<p\leq1$) which have very different behavior from the super-linear case. 
The super-linear  nonlinear eigen value problem can be described by the solutions of fast diffusion equations, where the solution will extinct at the finite time. So the Harnack type estimate plays an important role to analyze the asymptotic behavior  near the finite extinct time. 
On the other hand, the solution of sub-linear  eigenvalue problem will be approximated by the solutions of  slow diffusion equation, where the parabolic solution exists for all time. This difference allows us to have different method based on barriers and then have sharper results than the super-linear case. 
When $F$ is Laplace operator, the asymptotic behavior of the solution in the  degenerate or singular diffusion has been studied by many authors,
 Aronson , Berryman, Bonforte, Carrillo,  Friedman,  Galaktionov, Holland, Kamin,  Kwong, Peletier,  Toscani ,  Vazquez,  et al. We refer \cite{Va} for its details and references.

We also show that  the geometric property can be preserved in the degenerate fully nonlinear flow under the concavity condition of the operator and  hence such property also holds for the limit $\vp.$ 
To study the concavity of a solution,   the second difference of $u(x,t)$
$$C(x,y;u)=C(x,y)=2\left(u(x)+u(y)\right)-u\left(\frac{x+y}{2}\right) $$
is considered. 
Lastly, we show the eventual concavity of parabolic flow which means that the parabolic solution itself has such geometric property in finite time.

This analysis gives us sharp description of the asymptotic profile of the parabolic flow and affirmative answer for the well-known question on the convexity of level sets of the solution when the domain is convex.
We refer \cite{LV2}, \cite{KsL},\cite{GG}  for the detailed  history on the geometric issue. 

\subsection {}

Let $\mS^{n\times n}$ denote the set of $n\times n$ symmetric matrices and the norm  of a  matrix ,$||M||$, for $M\in\mS^{n\times n}$ is defined as the maximum absolute value among eigenvalues of $M.$ 
For  $0<\lambda \leq \Lambda$ (called ellipticity constants),  the Pucci's extremal operators, that play a crucial role in the study of fully nonlinear elliptic equations, are defined as ,for  $M\in\mS^{n\times n},$
\begin{align*}
\cM^+_{\lambda,\Lambda}(M)&=\cM^+(M)=\sup_{A\in\cA_{\lambda,\Lambda}}[tr(AM)], \\ \cM^-_{\lambda,\Lambda}(M)&=\cM^-(M)=\inf_{A\in\cA_{\lambda,\Lambda}}[tr(AM)],
\end{align*}
where $\cA_{\lambda,\Lambda}$ consists of the symmetric matrices  whose the eigenvalues  lie in $[\lambda,\Lambda]$.
We note that when $\lambda=\Lambda=1,$  the Pucci's extremal operators $\cM^\pm$ simply coincide with the Laplace operator.

In this paper, we assume that the nonlinear operator  $F : \mS^{n\times n} \rightarrow \re$
satisfies the following hypotheses  unless it is specifically mentioned :

$(F1)$ $F$ is a uniformly elliptic operator;  for all $M,N\in\mS^{n\times n},$
\begin{align*}
\cM^-(M-N) \leq F(M)-F(N) \leq\cM^+(M-N).
\end{align*}

$(F2)$ $F$ is positively homogeneous of order one; for all $t\geq 0$ and $M\in\mS^{n\times n}$
$$F(tM)=tF(M).$$

In addition, we assume that

$(F3)$ $F$ is concave.

The concavity condition of $F$ will be required when we show geometric property of parabolic flows. {The Pucci's extremal operator $\cM^-$ is one of nontrivial examples of the operator satisfying (F1), (F2) and (F3).} We may extend  $F$ on $\re^{n\times n}$ by defining $F(A):=F\left(\frac{A+A^T}{2}\right)$  for a nonsymmetric matrix $A.$

{Throughout this paper, we  assume that $\Omega$ is  a bounded domain with a smooth boundary in $\R^n  .$
}


We consider  viscosity solutions of \eqref{eq-intro-1},\eqref{eq-intro} which are proper notion of the weak solution for the fully nonlinear uniformly elliptic equation.   A continuous function $u \in\Omega$
 is said to be a \emph{viscosity subsolution} (respectively,
\emph{viscosity supersolution}) of $F(D^2 u(x))= f(x)$ in $\Omega$
 when the following condition holds: for any $x_0\in\Omega$ and $\phi\in C^2(\Omega)$ such that $u-\phi$ has a local maximum at $x_0$, we have
$$F(D^2\phi(x_0)) \geq f(x_0)$$
(respectively, if $u-\phi$ has a local minimum at $x_0$, we have
$F(D^2\phi(x_0)) \leq f(x_0)$). We say
that u is a viscosity solution of $F(D^2 u(x))= f(x)$ in $\Omega$ when it is both a viscosity subsolution and
supersolution. Viscosity solutions have been used to prove existence of solutions by Perron's method via the comparison principle.
We refer the details and regularity theory of the viscosity solutions to \cite{CIL},\cite{CC}.


\subsection{}\label{sec-intro-issues}
When the operator is fully nonlinear, there are several crucial issues to be discussed.
\begin{enumerate}[(i)]
\item Parabolic approach relies on the convergence of  the parabolic flows, $u(x,t)$,  to  eigen functions, $\vp(x),$ after normalizations, \eqref{eq-rn-1},\eqref{eq-rn-d}.
For nonlinear parabolic flow of divergence type, some key steps for the analysis of asymptotic behavior of  are based on the integration by parts, for example the existence of monotone integral quantities,\cite{Va}, which can not be applicable to the fully nonlinear operator.

On the other hand,  asymptotic Analysis in nondivergence form  can be  achieved in a couple of steps.
First, it is crucial to find an exact decay rate of $u(x,t)$, which will give us the right normalization of $u(x,t)$ so that the normalized parabolic flows converge to eigen functions, $\vp(x)$.
In fact, the exact decay rate is related to the first eigen value when $p=1$ and $m=1.$ We show the existence of the unique limit of normalized flow, $v(x,t)$, of $u(x,t),$ at Proposition \ref{lem-asymp}. 
When $0<p< 1$( or $1< m=\frac{1}{p}<\infty$), we prove Aronson-Benilan type estimate, Lemma \ref{lem-ab}, for degenerate fully nonlinear parabolic flows, which will give us almost monotonicity of $u(x,t).$ The uniqueness of the limit   of normalized flows is proven at Proposition \ref{prop-pme-app}.

\item Finally, we need to show that geometric properties of $u(x,t)$ will be preserved under the fully nonlinear parabolic flows, \eqref{eq-intro-1}. Geometric computation requires sophisticated computations to construct geometric quantities which satisfies maximum principle at Lemmas \ref{lem-c2f}, \ref{lem-c2f-p}.
The log-concavity of $u$ for $p=1$, the square-root concavity  of its pressure $u^{m-1}$ for $0<p<1$,  turn out to be preserved geometric quantities.
The difference of exponents comes from the difference in  homogeneities of the operators, \cite{Le}.

\end{enumerate}

\subsection{}
 This paper is organized into four parts as follows. 
At Section \ref{sec-nle}, we summarize the known facts about fully nonlinear uniformly
elliptic or parabolic partial differential equations.
And  at Section \ref{sec-nlev}, we show  
 some known results about the fully nonlinear elliptic eigenvalue problem $\eqref{eq-intro}$ and
the existence results of  positive eigen-functions
 for fully nonlinear elliptic problem as well as solutions
 of the nonlinear diffusion equations
in the range $0<p\leq 1.$

In Section \ref{sect-ufnl},  we deal with the fully nonlinear uniformly parabolic and elliptic equations and we discuss the log concavity of the first eigenfunction of nonlinear elliptic problem.
First, Bernstein's technique gives uniform
estimates for normalized solutions $v(x,t)= e^{\mu t} u(x,t)$
and then we use it to get the eigen-function as the limit of $v(\cdot,t)$ at Proposition \ref{lem-asymp}.
On the other hand,
we can choose initial data for this evolution having the desired
geometric property, and then
the evolution preserves the
geometric property.
Therefore the result for the elliptic problem will be
obtained in the limit $t\to\infty.$

Finally at Section \ref{sect-pme},  we show the long-time behavior of
the parabolic flow for $0<p<1$, Proposition \ref{prop-pme-app}. It is also proved that the pressure of the solution   preserves  square-root concavity under the parabolic flow and hence the concavity of eigenfunction is  proved.

$\textbf{Notations}$. Let us make a summary of the notations and definitions that  will be used .

$\bullet$ We denote by $\D u$ or $Du$ the spatial gradient of a function $u(x,t)$, and by $D^2u$ the Hessian
matrix. $D_ef$ denotes the directional derivative in the direction $e\in\mathbb{S}^{n-1}$.

$\bullet$ The expressions $D^2u\geq 0$, $D^2u\leq0$ are understood in the usual sense of quadratic forms.

$\bullet$ In order to avoid confusion between coordinates and partial derivatives, we will use the
standard subindex notation to denote the former, while partial derivatives will be denoted
in the form $f_{\al}$ for $\frac{\p f}{\p\al}=\p_{\al}f$. In general, $f_{\al}=\D_{e_{\al}}f$ for a unit direction $e_{\al}\in\mathbb{S}^{n-1}$
with a parameter $\al$. And second partial derivatives will be denoted
in the form $f_{\al\be}$ for $\frac{\p^2 f}{\p\al\p\be}=e_{\be}^TD^2fe_{\al}$. If the computation is invariant under the rotation, we may assume that
$\al=1,\cdots,n$ and that $\{e_1,\cdots,e_n\}$ is an orthonormal basis. This notation is usual in some
parts of the physics literature. But we will write $f_{\nu}$ and $f_{\tau}$ for the normal and tangential
derivatives since no confusion is expected.

$\bullet$ h.o.t. means 'higher order terms'.

\section{Preliminaries}\label{sec-nle}
For the reader's convenience, we are going to summarize basic facts and estimates for elliptic
fully nonlinear equation  $F(D^2u)=f(x)$ in a bounded domain
$\Omega\subset \R^n$, \cite{CC, CIL} and for parabolic fully nonlinear
equations $u_t=F(D^2u)+f(x)$ in $Q_T=\Omega\times(0,T]$, \cite{CIL, L,
  W1, W2}, 
where $F$ satisfies the condition (F1).
\begin{enumerate}
\item The {\it existence and  uniqueness of the viscosity solution}, 
  {\it comparison principle} between super- and sub-solutions, {\it minimum
  principle and maximum principle} in 
  elliptic or parabolic Dirichlet problem, and their references can be found at \cite{CIL,
    CC, W1}.
\item The {\it strong maximum principle} holds for $f=0$ at elliptic equation,
  Proposition 4.9, \cite{CC} and the same argument with Corollary
  3.21, \cite{W1}, gives us the strong maximum principle for the
  parabolic equation. A version of  strong maximum principle for fully
  nonlinear equation with nonhomogeneous operator  has been
  proved at Lemma \ref{lem-cp}. The strong maximum principle for elliptic or
  parabolic equation says that
  whenever a subsolution $u$ touches a super-solution $v$ from below at an
  interior point , $u\equiv v$ on the domain $\Omega$ or $Q_T,$
  respectively.
\item{[{\it Local Regularity}]} We
  refer the regularity theory for elliptic equation to \cite{GT, CC}
  and parabolic case to \cite{L,W1,W2}.
  In \cite{CC}, we can find 
H\"older
  continuity ($k=0,0<\alpha<1$)  at Proposition 4.10, 
  Local $C^{1,\alpha}$-regularity($k=1,0<\alpha<1$)  at  Theorem 8.3, Local $C^{1,1}$-regularity($k=2,\alpha=0$) for convex or
  concave operator $F$ at Proposition 9.3, local
  $C^{2,\alpha}$-regularity($k=2,0<\alpha<1$) for H\"older continuous
  $f$  at Theorem 8.1, local $C^{\infty}$-regularity  for smooth $F$ and
  $f$. When $F$ and $f$ is analytic, $u$ will be analytic following
  Theorem 10 at Section 2.2, \cite{E}.
\item{[{\it Global
      Regularity}]} We also refer the global Schauder theory ,$$\|u\|_{C^{k+2,\alpha}(\overline{\Omega})}\leq
  C(\|u\|_{L^{\infty}(\overline{\Omega}
    )}+\|f\|_{C^{k,\alpha}(\overline{\Omega})})$$ to Theorem 9.5 ,
  \cite{CC}.  Therefore if $\partial\Omega$ is $C^{2,\alpha}$-surface,
  then the viscosity solution will be classical. The similar results
  hold for parabolic equation, \cite{L}.
\item{[{\it Harnack Inequality}]} {\it The Harnack inequality} for
  a nonnegative elliptic solution is the following, Theorem 4.8, \cite{CC}:  for a nonnegative elliptic  solution $u$ in $B_3$,
  we have $\displaystyle\sup_{B_1} u\leq C(\inf_{B_1}u+\|f\|_{L^n(B_2)}
)$
for  a uniform constant $C>0$. Similar parabolic
version can be found at \cite{W1}.
\end{enumerate}

\section{Nonlinear eigenvalue problem}\label{sec-nlev}

In this section we are going to study solutions to the fully nonlinear elliptic eigenvalue problem
\begin{equation}
\begin{cases}
 F(D^2 \phi(x)) = -\mu\phi^p(x)\quad\text{in $\Omega,$}\\
  \phi >0\quad\text{in $\Omega,$}\\
 \phi = 0\quad\text{on $\partial\Omega,$}
\end{cases}
\label{NLEV0}
\tag{{\bf NLEV}}
\end{equation}
where $\Omega$ is a smooth bounded domain in $\re^n$ and
$F$ is a uniformly elliptic and positively homogeneous operator of order one  defined on $\mS^{n\times n}.$
  First, let us introduce the existence theorem of the  positive eigen-function that was proven by Ishii and Yoshimura.  The simplified proof can be found at \cite{A}.

\begin{theorem}\cite{IY}
Suppose that $F$ satisfies $(F1)$ and $(F2)$ and that $\Omega$ is a smooth bounded  domain in $\re^n.$
Then there exist {$\vp \in C^{1,\al}(\overline\Omega)  , (0<\alpha<1)$ } and $\mu>0$ such that  $\vp>0$ in
$\Omega$ and $\vp$ satisfies
\begin{equation}
\begin{cases}
 -F(D^2 \vp(x)) =\mu \vp(x)\quad\text{in $\Omega,$}\\
\vp(x)=0\quad\text{on $\partial\Omega.$}
\end{cases}
\label{eq-main-1}
\tag{{\bf EV}}
\end{equation}
Moreover, $\mu$ is unique in the sense that if $\rho$ is another
eigen-value of $F$ in $\Omega$ associated with a nonnegative
eigen-function, then $\mu=\rho$ ; and is simple in the sense of that
if $\psi$ in $C^0(\overline\Omega)$ is a solution of \eqref{eq-main-1} with
$\psi$ in place of $\vp$, then $\psi=c\vp$ for some $c \in \re$.
\label{thm1}
\end{theorem}



Now  let us state the Hopf's Lemma that will be used frequently  when we compare a solution with barrier.

\begin{theorem}[Hopf's Lemma]\label{hopf} 
Suppose that $\Omega$ satisfies an interior sphere condition. 
Let $u\in C(\overline{\Omega})$ be a nonzero 
 viscosity supersolution of
$$\cM^-(D^2 u)\leq 0\,\,\,\,\mbox{in}\,\,\Omega.$$ 
Then for $x_o\in\p\Omega$ such that $u(x)>u(x_o)$ for all $x\in\Omega,$ we have 
$$\displaystyle\liminf_{x\in\Omega\to x_o}\frac{u(x) -u(x_o)}{|x-x_o|} > 0.$$
Especially, if the outer normal derivatives of $u$ at $x_o$ exists, then
$$\displaystyle\frac{\p u}{\p \nu}(x_o)<0,$$
where $\nu$ is the outer normal  vector to $\p\Omega$ at $x_o.$\\
In particular, if $u=0$ on $\p\Omega,$ $u>0$ in $\Omega$ and $u\geq M>0$ on $\Omega'\Subset\Omega,$ then
$$\displaystyle\frac{\p u}{\p \nu}(x_o)<- c_o(M,\,\mbox{dist}\,(\Omega',\Omega)).$$
\end{theorem}
 We refer to Lemma 3.4 at \cite{GT} for  uniformly elliptic linear equation and
  Lemma 2.6 at \cite{L} and Appendix at \cite{A} for uniformly parabolic fully nonlinear equation. The Hopf's
  lemma for uniformly elliptic fully nonlinear equation follows by the
  comparison between super-solution and a barrier
  $R^{-\alpha}-|x|^{-\alpha}$ for large $\alpha>0$ and a small $R>0$ as  Lemma 2.6, \cite{L}. 
Hopf's Lemma  for the parabolic equation \cite{L}  holds in the following way :
$$\displaystyle\liminf_{x\to x_o, s\to t}\frac{u(x,s) -u(x_o,t)}{\sqrt{|x-x_o|^2+(t-s)}} > 0$$
for any $x\in\Omega$ and $s\leq t.$ 
 
\subsection{Case $0< p < 1$}

In this case, we consider the following equation
\begin{equation}
\begin{cases}
 -F(D^2 \ f^m(x)) = \frac{1}{m-1}f(x)\quad\text{in $\Omega$},\,\, m>1,\\
 f = 0\quad\text{on $\partial\Omega,$}\\
 f > 0\quad\text{in $\Omega,$}
\end{cases}
\label{eq-main-1-pme}
\end{equation}
which is the asymptotic profiles of the equation
\begin{equation}\label{eqn-gen}
\left\{
\begin{array}{ll}
H[u]=u_t(x,t)-F(D^2 u^m(x,t))=0\qquad &\text{in $Q=\Omega\times (0,T]$},\\
u(x,0)=u_o(x), & \text{in $ \Omega $},\\
u(x,t)=0\quad &\text{on $\p\Omega\times (0,T].$}
\end{array}\right. \qquad \qquad
\end{equation}
We assume  $u_o$ has nontrivial bounded  gradient on $\p\Omega,$ i.e, 
 $$u_o^m\in\cC(\overline\Omega),$$
where
$$\cC(\overline\Omega):=\{ h\in C^o(\overline\Omega) |
c_o\,\mbox{dist}(x,\p\Omega) \leq h(x) \leq
C_o\,\mbox{dist}(x,\p\Omega)\,\mbox{for } 0<c_o\leq C_o<+\infty\}.$$
If we set $\vp = f^m$ and $p=\frac{1}{m},$ then $\vp$ is the solution of $\eqref{NLEV0}$ with an eigenvalue $\frac{1}{m-1}.$
For the sub-linear case, $0<p<1,$ we have the following comparison principle and the existence and uniqueness result of nonlinear eigenfunction.


\begin{lemma}[ Comparison principle]\label{lem-cp}
Suppose $F$ satisfies $(F1), F(0)=0$ and that  either $(F2)$  or $(F3).$
Let $v$ and $ w$ be in $ C^2(\Omega)\cap{ C^1(\overline{\Omega})}$ such that $v,w\geq 0.$
If  $F(D^2v)+\frac{1}{m-1}v^{\frac{1}{m}}\leq 0 \leq F(D^2w)+\frac{1}{m-1}w^{\frac{1}{m}}$ in $\Omega$
and if $v \geq w\,$ on $\partial\Omega,\,\,$
 then  $v \geq w\,$ in $\Omega.$
\end{lemma}

\noindent {\sc Proof.}
Suppose that $v < w$ for some point in $\Omega.$ Since $v$ satisfies
$$\cM^-(D^2v) \leq F(D^2v) \leq -\frac{1}{m-1}v^{\frac{1}{m}}  \leq0,\,\,$$
we have  $v > 0\,$ in $\Omega\,$ and $|\D v|> 0\,$ on $\partial\Omega$ by the strong minimum principle and Hopf's lemma \ref{hopf}.
Let $t^* = \,\mbox{inf} \{ t>0 | v < tw \,\,\mbox{for some point in}\,\,\Omega\}.$ Then $0< t^* <1.$
Set $z = v - t^*w ,$ and then the nonnegative function $z$ vanishes at
  some point in $\overline{\Omega}$  
and $z$ satisfies 
\begin{align*}
\cM^-(D^2z) \leq F(D^2v)-F(D^2t^*w) 
& \leq \frac{1}{m-1}(t^*w^{\frac{1}{m}} - v^{\frac{1}{m}})\\
&\leq \frac{1}{m-1}((t^*w)^{\frac{1}{m}} - v^{\frac{1}{m}}) \leq 0.
\end{align*}
Assume that $z\not\equiv 0.$
From the strong minimum principle and Hopf's lemma,
 we have  $z > 0\,$ in $\Omega\,$ and $|\D z|> 0\,$ on $\partial\Omega.\,\,$
Then we can choose $\e > 0 \,$ such that {$z - \e v \geq 0\,\,$ in $\Omega.$}
It's a contradiction to the definition of $t^*.$ Thus we get  $z\equiv 0$ and  $v=t^*w$ in $\overline\Omega.$

(i)  First, let us assume that $v$ is a strictly supersolution, i.e., $F(D^2v)-\frac{1}{m-1}v^{\frac{1}{m}}<0,$
 we have 
\begin{align*}
0 > & F(D^2v)+\frac{1}{m-1}v^{\frac{1}{m}}\geq t^*F(D^2w)+\frac{1}{m-1}(t^*w)^{\frac{1}{m}}\\
&= t^*\left\{F(D^2w)+\frac{t{^*}^{\frac{1}{m}-1}}{m-1}w^{\frac{1}{m}}\right\}  \geq t^*\left\{F(D^2w)+\frac{1}{m-1}w^{\frac{1}{m}}\right\} \geq0,
\end{align*}
which is a contradiction.

(ii) Now,  assume that $v$ is a supersolution, i.e., $F(D^2 v)\leq -\frac{1}{m-1}v^{\frac{1}{m}}.$ Then, we have that  $v>0$ in $\Omega$ by the strong minimum principle and Hopf's lemma.

Let $v^{\e} := (1+\e)v$   for $\e>0.$  Then 
$v^{\e}$ satisfies 
\begin{align*}
F(D^2v^{\e})+ \frac{1}{m-1}(v^{\e})^{\frac{1}{m}}&\leq(1+\e)F(D^2v)+ \frac{(1+\e)^{\frac{1}{m}}}{m-1}v^{\frac{1}{m}}\\
&\leq \frac{1}{m-1}v^{\frac{1}{m}}\{(1+\e)^{\frac{1}{m}}-(1+\e)\} <0,
\end{align*}
i.e., $v^{\e}$ is  a strictly supersolution.
By (i) , we get  $v^{\e}=(1+\e)v\geq w$ in $\Omega.$ Letting $\e\to 0,$ we have 
$v\geq w$ in $\Omega.$
\qed

\begin{theorem}\label{thm-phi-m}
Suppose $F$ satisfies $(F1)$ and $(F2)$. The nonlinear eigenvalue problem has a unique positive viscosity solution $\phi\in C^{0,1}(\overline{\Omega})\cap C^{1,\al}({\Omega}),$
i.e.,
\begin{equation}
\begin{cases}
 -F(D^2 \phi(x)) = \frac{1}{m-1}\phi^{p}(x)\quad\text{in $\Omega,$}\\
 \phi = 0\quad\text{on $\partial\Omega,$}\\
 \phi > 0\quad\text{in $\Omega,$}
\end{cases}
\label{NLEV}
\tag{{\bf NLEV}}
\end{equation}
where $p= \frac{1}{m}$. The eigen-function $\phi$ satisfies $\displaystyle\inf_{\p\Omega}|\D\phi|\geq\delta_o>0.$
Moreover, if $F$ is $C^1$, $\phi$ is of $ C^{\infty}({\Omega}).$
\label{thm2}
\end{theorem}

\noindent {\sc Proof.}
(i) The uniqueness of the solution follows from Comparison Principle.
It suffices to establish the existence of positive super and sub-solutions with zero boundary value in order to prove the existence of the solution.
Let $h$ be the solution of
\begin{equation}
\begin{cases}
 F(D^2 h(x)) = -1\quad\text{in $\Omega,$}\\
 h = 0\quad\text{on $\partial\Omega,$}\\
 h > 0\quad\text{in $\Omega.$}
\end{cases}
\end{equation}
If  we select $t > 0\,$  satisfying $t^{1-\frac{1}{m}}||h||_{L^{\infty}(\Omega)}^{-\frac{1}{m}} = \frac{1}{m-1},$ then
\begin{align*}
F(D^2(t h)) = -t^{1-\frac{1}{m}}h^{-\frac{1}{m}}(t h)^{\frac{1}{m}} \leq  -\frac{1}{m-1}(t h)^{\frac{1}{m}}
\end{align*}
i.e., $h^+ :=t h$ is a super-solution.

On the other hand, let $\vp$ be the first  eigen-function of $\eqref{eq-main-1}.$ Choose $s>0$ so that $\mu(s ||\vp||_{L^{\infty}(\Omega)})^{1-\frac{1}{m}} \leq \frac{1}{m-1}$, then
$F(D^2(s\vp)) 
\geq -\frac{1}{m-1}(s\vp)^{\frac{1}{m}}.$
Thus $h^- := s\vp $ is a sub-solution.

Thus the comparison principle Lemma \ref{lem-cp} gives that $h^-\leq h^+$ and there is a viscosity solution $\phi$ such that 
$h^-\leq \phi \leq h^+$ from \cite{CIL}. Since $\phi^{p}\in L^{\infty}(\Omega),$ {$\phi$ is  of $C^{1,\alpha}(\Omega)$} from the regularity theory, \cite{CC}. Since $ F(D^2 \phi(x)) = -\frac{1}{m-1}\phi^{p}(x)\leq 0,$ $\phi$ satisfies $\displaystyle\inf_{\p\Omega}|\D \phi|\geq\delta_o>0$ from Hopf's lemma.

(ii) 
Now we are going to show $\phi\in C^{0,1}(\overline{\Omega})\cap C^{1,\alpha}(\Omega).$ 
 First, we note that
there are $0<c_o\leq C_o<\infty$ such that $c_o\,\mbox{dist}\,(x,\p\Omega)\leq h^-\leq \phi\leq h^+\leq C_o\,\mbox{dist}\,(x,\p\Omega)$ from Hopf's Lemma for $h^-$ and $C^{0,1}(\overline{\Omega})$ - regularity of $h^+,$ \cite{CC}.

Let $\delta>0$ be a constant such that $B_{\delta}(x)\subset\Omega$ for $\mbox{dist}\,\,(x,\p\Omega)>\delta.$ 
 For  $x_o\in\Omega$ such that $\mbox{dist}\,\,(x_o,\p\Omega)<\delta,$ set $\,\mbox{dist}\,(x_o,\p\Omega)=2\e.$ Now we scale the function $\phi,$ 
$$\phi_\e(x)=\frac{1}{\e}\phi(x_o+\e x).$$
Then $0<c_o\leq \phi_\e(x)\leq 3C_o$ in $B_1(0)$ and $\phi_\e$ satisfies 
$$F(D^2\phi_\e)=-\frac{\e^{1+p}}{m-1}\phi_\e ^p \in L^{\infty}(B_1(0))\quad\mbox{uniformly}.$$
From the regularity theory,\cite{CC}, we have 
$$|D\phi(x_o)|=|D\phi_\e(0)| \leq \tilde C\,\,\mbox{for some uniform constant }\,\,\tilde C>0,$$
Therefore, we have $|D\phi(x_o)| \leq \tilde C $ and we  deduce that $\phi\in C^{0,1}(\overline\Omega).$

(iii) When $F$ is $C^1, $ the operator becomes a  linear operator from the positive homogeneity  of order one. Thus the result follows.        
\qed

We  state  the following comparison principle  of the  solution $,u\,,$ of the parabolic flow $\eqref{eqn-gen}$ for the case $m>1$ and 
we consider the following equation:
\begin{equation}\label{eqn-gen1}
\left\{
\begin{array}{ll}
F(D^2 v(x,t))=(v^{\frac{1}{m}})_t(x,t) \qquad &\text{in $Q_T=\Omega\times (0,T]$},\,\,\,m > 1,\\
v(x,0)=v_o(x)=u_o^{{m}}\in \cC(\overline\Omega), & \\
v(x,t)=0\quad &\text{on $x\in \partial\Omega$}.
\end{array}\right. \qquad \qquad
\end{equation}
The proof of Comparison principle for the case $m>1$, is  the same as the case $\textsf{m}_{\Omega,F} <m<1$, \cite{KsL}.
The similar argument as  Lemma \ref{lem-cp} gives us the following Lemma.

\begin{lemma}[ Comparison principle]
Suppose $F$ satisfies $(F1), F(0)=0$ and either $(F2)$ or $(F3).$ Let  $v\,,w \in C^{2,1}(Q_T)\cap C^0(\overline{Q_T})\,$ such that $v,w >0$ in ${Q_T}.$  If $F(D^2 v)-(v^{\frac{1}{m}})_t  \leq 0 \leq F(D^2 w)-(w^{\frac{1}{m}})_t $ in $Q_T$
  and if $v \geq w\,$ on $\p_p Q_T=\left(\Omega\times\{0\} \right) \cup\left( \p\Omega\times(0,T]\right),$
 then  $v \geq w\,$ in $Q_T.$
\end{lemma}




\begin{theorem}\label{thm-ex-pme} Suppose $F$ satisfies $(F1), F(0)=0$ and either $(F2)$ or $(F3).$ Let $u_o^m$ be in $\cC(\overline\Omega).$ When $m>1,$ there exists a unique solution of porous medium type \eqref{eqn-gen}. Moreover, the solution $u$ is positive in $\Omega\times(0,+\infty).$
\end{theorem}
\noindent {\sc Proof.}
 Let $f=\phi^{\frac{1}{m}}$ for $\phi$ in Theorem \ref{thm-phi-m}.  First, we note that $0$ and $f(x)(k+t)^{-\frac{1}{m-1}}\,\,$ (for any $k>0$) are solutions of
$u_t(x,t)=F(D^2 u^m(x,t))\,$ in $\,
Q_T,$ with zero boundary. 

We  construct a supersolution using self-similar solutions. Let   $\phi^+$ be an eigen-function with the Pucci's operator $\cM^+$ in Theorem \ref{thm-phi-m}.
For a given $\e>0,$ we can choose $K>0$ such that $0 < u_o^m(x) \leq \phi^+(x)K^{-\frac{1}{m-1}}$ since  $\displaystyle \inf_{\p\Omega}|\D \phi|>0.$ 
Then ${\phi^+}^{\frac{1}{m}}(x)(K+t)^{-\frac{1}{m-1}}\,\,$ is a supersolution of \eqref{eqn-gen} with any  $F$ as the  operator since $\cM^-\leq F \leq \cM^+.$ Therefore  
there exists a unique
solution $u\,$ of $\eqref{eqn-gen}.\,\,$ Moreover, $u$ satisfies $$0\leq u(x,t) \leq {\phi^+}^{\frac{1}{m}}(x)(K+t)^{-\frac{1}{m-1}}\,\,$$
in $\,Q_T$ from the Comparison principle. 

{In addition, we are going to show that $u>0$ in $Q_T$ if $u_o^m\in\cC(\overline{\Omega}).$  
Let $\Omega'$ be any smooth compact subset of $\Omega$. From Theorem
\ref{thm-phi-m}, there is a positive eigenfunction $\vp_1$
corresponding to $\Omega'$ with the operator $\cM^-.$   Set $g(x)=\vp_1^{1/m}$ and then
$U_1=g(x)(K+t)^{-\frac{1}{m-1}}\,\,$ solves 
\eqref{eqn-gen} in $\Omega'\times(0,T]$ with the operator $\cM^-.$ Since $u_0(x)>0$ on a compact set $\overline{\Omega}'$
and $u_0(x)$ is continuous, there is  a large $K>0$ such that
$U_1(x,0)=g(x)K^{-\frac{1}{m-1}}\leq u_0(x)\,\,$ on
$\overline{\Omega}'$. From the comparison principle in $\Omega'$, we
have $u(x,t)\geq U_1(x,t)>0$ on $\Omega'\times(0,T]$. By taking $\Omega'$ arbitrary, we have $u>0$ in $\Omega\times[0,\infty)$.}
 \qed


\section{Uniformly fully nonlinear equation   }\label{sect-ufnl}

We consider the solutions $u(x,t)$ of the problem
\begin{equation}\label{eqn}
\left\{
\begin{array}{ll}
H[u]=u_t(x,t)-F(D^2 u(x,t))=0\qquad &\text{in $Q=\Omega\times (0,+\infty)$},\\
u(x,0)=u_o(x)\in C^o(\overline\Omega), & \\
u(x,t)=0\quad &\text{for $x\in \partial\Omega\times(0,+\infty)$},
\end{array}\right. \qquad \qquad
\end{equation}
where $\Omega$ is a bounded domain of $ \re^n$ {with a  smooth
boundary.}

\subsection{Asymptotic Behavior}
In this subsection, we are going to analyze the asymptotic behavior of the solution $u$ of \eqref{eqn}. First, we will find the exact decay rate of $u$  comparing it with  barriers constructed by using the principal eigen-value, $\mu$,  and a positive eigen-function, $\vp(x)$ .
\begin{lemma}
Suppose $F$ satisfies $(F1)$ and $(F2)$.
For any positive $u_0\in \cC(\overline\Omega),$ 
there are  $0<C_1\leq C_2$   such that
\begin{equation*}
C_1\vp(x)e^{-\mu t}<u(x,t)<C_2\vp(x)e^{-\mu t},
\end{equation*}
 for $t>0$.
\end{lemma}

\noindent {\sc Proof.}
By Hopf's lemma and { $C^{0,1}- $ regularity of $\vp$,} 
 we have $0<|\D\vp| <+\infty\,\,\,\text{on $\partial\Omega$}.$ So we can choose $C_2 > C_1>0$ such that $C_1\vp(x)<u_o(x)<C_2\vp(x)\quad\text{in $\Omega$}.$
Since $C\vp(x)e^{-\mu t}$ is a solution of \eqref{eqn} for any constant $C>0$, the comparison principle gives us the result.
\qed

\begin{lemma}
\label{lem-asym} Suppose that $F$ satisfies $(F1)$ and $(F2).$ 
For any nonnegative and nonzero $u_o\in C^0(\overline\Omega)$, there is  $t_0>0$ such
that $$C_1\vp(x)<u(x,t_0)<C_2\vp(x),$$
 for some $0<C_1\leq C_2$
and then  for $t\geq t_0$
$$C_1\vp(x)e^{-\mu t}<u(x,t)<C_2\vp(x)e^{-\mu t}.$$
\end{lemma}

\noindent {\sc Proof.}
We are going to construct a subsolution of $\eqref{eqn}$ which expands in time.
Define $g(x,t)\,=\,\frac{1}{t^{\be}}\exp\left(-\al\frac{r^2}{t}\right)$, where $ \al =\frac{1}{4\lambda},\,\,\be = \frac{\Lambda n}{2\lambda}\,\,\mbox{and}\,\, r = |x|.$
We can easily see that at the point $(r,0,\cdots,0),$
\begin{align*}
\p_{ij}g &= 0\qquad\mbox{if}\,\,\,i\neq j,\\
\p_{11}g &=2\al \frac{g}{t^2}(2\al r^2-t),\\
\mbox{and}\,\,\,\p_{ii}g &=-2\al \frac{g}{t}\qquad\mbox{if}\,\,\,i>1.
\end{align*}
Then we check for $r^2<\frac{t}{2\al}$
\begin{align*}
\cM^- (D^2g) - g_t = \Lambda \p_{11}g+(n-1)\Lambda\p_{22}g-g_t= \frac{g}{t^2}\{ t(\be-2\al\Lambda n)+\al r^2(4\Lambda\al-1)\}\geq 0
\end{align*}
and for $r^2\geq\frac{t}{2\al}$
\begin{align*}
\cM^- (D^2g) - g_t &=\frac{g}{t^2}\{t[\be-2\al(\lambda+(n-1)\Lambda]+\al r^2(4\lambda\al-1)\}\geq 0.
\end{align*}
Thus $g$ is a subsolution of $F(D^2u)-u_t=0.$  For positive constants $\tau_o,\,c_o\,$ and $\delta_o,$
we define 
$$h(x,t) := \mbox{max}\left\{\displaystyle c_o\frac{1}{(t+\tau_o)^{\be}}\exp\left(-\al\frac{|x-x_0|^2}{t+\tau_o}\right)-\delta_o\,,\,0\right\}$$
and then $h$ is also a subsolution as long as $\mbox{supp}h(\cdot,t) \subset \Omega.$

Since $u_o \not\equiv 0,\,$ there exists a point $x_o \in\Omega$ such that $u(x_o) := m_1> 0.$ We choose $\rho > 0 $ and $\eta>0$ small
so that $B_{\rho}(x_o) \subset\subset \Omega,\,\,  u_o(x)\geq \frac{m_1}{2}\,= m_o \,\mbox{in}\,\, B_{\rho}(x_o),\,$
and $2\rho  < \mbox{dist} (x_o , \partial \Omega)$ 
 and   that $0<\eta \leq \rho ,$ and $B_{2\eta}(y) \subset \Omega\,$ for $y \in \{x \in\Omega\,|\, \mbox{dist}(x,\partial\Omega) \geq 2\eta \} \equiv \Omega_{2\eta}.$
By taking $\tau_0, c_0$ and $\delta_0$ such that
\begin{align*}
\eta^2 = 4\Lambda n\tau_0 (> 2\Lambda n\tau_0),\quad\frac{c_0}{{\tau_0}^{\be}} \exp(-\al\frac{\eta^2}{\tau_0}) = \delta_0\quad\mbox{and}\quad\frac{c_0}{{\tau_0}^{\be}} - \delta_0 = m_0,
\end{align*}
then the support of $h(x,t)$ is increasing for $0<t \leq\displaystyle \frac{1}{e}\left(\frac{c_0}{\delta_0}\right)^{1/\be}-\tau_0$ with $h(x,0) < u_o(x).$
In fact, at time $t_0=\displaystyle\frac{1}{e} \left(\frac{c_0}{\delta_0}\right)^{1/\be}-\tau_0= \displaystyle\frac{e-1}{4\Lambda\lambda}\eta^2$, the support of $h(x,t)$  becomes the ball with radius $\displaystyle\sqrt{\frac{e}{2}}\,\eta$ centered at $x_o.$
Comparison principle implies $h(x,t) \leq u(x,t)\,\, \mbox{in} \,\,\Omega\times(0,t_0]$
and hence $u(x,t_0) > 0 $ in $B_{\sqrt{\frac{e}{2}}\,\eta}(x_o)$ at $ t_0=\displaystyle\frac{e-1}{4\Lambda\lambda}\eta^2 > 0.$

For any point $y\in\p\Omega,$ we have a chain of uniform number of balls with radius $\sqrt{\frac{e}{2}}\,\eta$ from $x_o$ to $y$
 and each ball will be filled by the above
subsolution $h$ starting at the previous ball. Since all of argument can be carried out at finite step
only depending on the initial data and the domain, there is a time $t_1$ such that
$u(\cdot,t_1) > 0 $ in $\Omega$ and $|\nabla u(y,t_1)| > 0$ for $y\in\p\Omega.$
Thus, there is $C_1>0$ such that $C_1\vp(x)e^{-\mu t_1}<u(x,t_1)$ in $\Omega.$ Since
 { $u$ is $C^{1,\beta}(\overline{\Omega\times[t_1,t_1+1]}),$ }   there is $C_2>0$ such that $u(x,t)<C_2\vp(x)e^{-\mu t}$ in $\Omega\times[t_1,\infty).$ Therefore, the result follows.
\qed

To refine the asymptotic behavior, let us introduce the normalized function 
\begin{equation}\label{eq-rn-1}
v(x,t)=e^{\mu t}u(x,t).
\end{equation} Then, $v(x,t)$ satisfies $v_t=F(D^2 v)+\mu v$ if the operator $F$ satisfies the condition $(F2)$ and we deduce the following Corollary from Lemma \ref{lem-asym}.

\begin{corollary}\label{cor-bd} Under the same assumption of  Lemma \ref{lem-asym},
$v(x,t)=e^{\mu t}u(x,t)$ has the following estimate: 
\begin{equation*} 
||v(x,t)||_{L^{\infty}(\Omega\times [t_0,\infty))}\leq C||{v(x,t_0)}||_{L^\infty(\Omega)},
\end{equation*}
where $t_0>0$ is in Lemma \ref{lem-asym}.
\end{corollary}

Before studying  fine asymptotic behavior of  parabolic solutions, let us  summarize the regularity theory of uniformly parabolic equation.
\begin{theorem} [{Global Regularity for $m=1$}]
Suppose that the domain $\Omega$ is bounded and smooth and $F$ satisfies {$(F1).$}
\begin{enumerate}[(i)]
\item  Let $u$ be a solution of $\eqref{eqn}$ and let $Q= \Omega\times(\delta_o,T)$ for any  
 $T>\delta_o>0.$
 \begin{enumerate}[(a)]
\item  $u$ is of $C^{1,\beta}(\overline Q)$ for some $0<\beta<1.$ 
\item  If $F$ is concave,  $u$ is of $C^{1,1}(\overline Q).$ 
\item  If $u\in C^{1,1}(\overline Q)$ and if $F$ is concave or convex,   $u$ is of $C^{2,\beta}(\ Q)$ for some $0<\beta<1.$
\item  If $u\in C^{2,\beta}(\overline Q)$ and  $F\in C^{\infty},$  $u$ is of $C^{\infty}(\overline Q).$
\end{enumerate}
\item Let $v(x,t)$ be  a bounded solution of $v_t=F(D^2 v)+\mu v.$ (a),(b),(c) and (d) for $v$ also hold. 
\end{enumerate}
\end{theorem}\label{thm-reg-1}
We refer the regularity theory    to  \cite{GT, CC, L,W1,W2}. We note that in this parabolic setting,  $C^{\beta}$ means that $C^{\beta}$ in $x$ and  $C^{\beta/2}$ in $t$.
 
Let us prove the  interior $C^{1,1}_x$ - estimate for reader's convenience through Bernstein's  computation.

\begin{lemma}
\label{lem-ck} Suppose that $F$ satisfies $(F1),$ $(F3)$ and $F(0)=0$ and { $F$ is of  $C^2.$ } Then, a bounded solution ${ v\in C^4}$ of $v_t=F(D^2 v)+\mu v ,(\mu\in\R)$  satisfies
\begin{align*}\label{eq-ck}
||v(x,t)||_{C^{1,\alpha}(Q_{\frac{1}{2}})}\leq C||v||_{L^\infty(Q_1)}\\
\text{and}\,\,||D^2v||_{L^{\infty}(Q_{\frac{1}{2}})}+{||v_t||}_{L^{\infty}(Q_{\frac{1}{2}})}\leq C||v||_{L^\infty(Q_1)}
\end{align*}where $Q_{R}:= B_R(0)\times(-R^2,0).$
Moreover if $F$ is smooth, 
$$||v(x,t)||_{C^k(Q_{\frac{1}{2}})}\leq C||v(x,t)||_{L^\infty(Q_1)}
$$
 for $k=1,2,\cdots.$
\end{lemma}
\noindent {\sc Proof.} (i) 
Let $M := ||v(x,t)||_{L^\infty(Q_1)}$ and let $\psi\in
C^{\infty}(\overline{Q_1})\,$ be a parabolic cutoff function such that $0\leq\psi\leq 1
\,\mbox{in}\,\,\overline{Q_1},\,
\psi=1\,\mbox{in}\,\overline{Q_{1/2}},\,\psi=0\,\mbox{on}\,\p_pQ_1\,\,$
and $|\psi|+|\D\psi|+|D^2\psi| < c=c(\psi).$

For large $\delta > 0 $ (to be chosen later), define
$$h(x,t)=\delta(M-v)^2 + \psi^2|\D v|^2+M^2\frac{8\delta|\mu|}{\lambda}\frac{x_1^2}
{2}.$$ Now, we define the uniformly elliptic operator
$$L[w]:=F_{ij}(D^2v)D_{ij}w,$$
and  the  uniformly parabolic operator $H[w] := L[w] - w_t$ and  we have that  { $H[v]\leq -\mu v, H[v_{e}]=-\mu v_e,$} and $H[v_{ee}]\geq -\mu v_{ee}$ from the concavity of $F$ and $F(0)=0$ using the function $a(s)=F((1-s)D^2 v)+(1-s)(\mu v-vt)$  as in the chapter 9 at \cite{CC}. Using Bernstein's technique, we  get
\begin{align*}
H[h] = Lh - h_t
&\geq 2\delta \lambda|\D v|^2 +2\delta(M-v)(-F_{ij}D_{ij}v+v_t)
+2|\D v|^2\lambda|\D\psi|^2  \\
&\quad-2\psi|\D v|^2|DF||D^2\psi| + 8\psi F_{ij}D_i\psi
D_kvD_{kj}v+2\psi^2\lambda|D^2v|^2\\
&\quad+2\psi^2D_kv(F_{ij}D_{kij}v -v_{kt})
-2\psi|\psi_t||\D v|^2+ M^2\frac{8\delta|\mu|}{\lambda}F_{11}\\
&\geq 2\delta \lambda|\D v|^2 +2\delta(M-v)\mu v +8M^2\delta|\mu|
+2|\D v|^2\lambda|\D\psi|^2  \\
&\quad-2\psi|\D v|^2|DF||D^2\psi| + 8\psi F_{ij}D_i\psi D_kvD_{kj}v
+2\psi^2\lambda|D^2v|^2+2\psi^2|\D v|^2\\
&\quad-2c|\D v|^2\geq 0
\qquad\mbox{{for large}}\,\,\,\delta=\delta (c(\psi), \Lambda, 
\lambda, n)>0.
\end{align*}
Since  $$\,\, h \leq \delta M^2 +
M^2\frac{8\delta|\mu|}{\lambda} \leq C M^2\quad\mbox{on}\,\,\,\p_p\Omega,$$ we obtain  that  $\displaystyle\sup_{Q_1}h \leq C M^2$ from the maximum
principle and  hence
$$||\D v(x,t)||_{L^{\infty}(Q_{1/2})}\leq C||v(x,t)||_{L^\infty(Q_1)}.$$

(ii) 
$||D^2v||_{L^{\infty}(Q_{\frac{1}{2}})}\leq C||v||_{L^\infty(Q_1)}$ comes from applying the maximum principle on
$$g=\delta (v_{e})^2+ \psi^2(v_{ee})^2+{\delta CM^2|\mu|x_1^2}$$
for any direction $e\in S^{n-1},$ as Proposition 9.3, \cite{CC}.
\qed

 Now, we are going to show  normalized parabolic flow $v(x,t)=e^{\mu t}u(x,t)$ has the unique limit as $t\ra \infty$  and     use the approach presented at \cite{AT} to obtan the uniqueness
 of the limit. 
\begin{proposition}\label{lem-asymp} Suppose $F$ satisfies $(F1)$ and $(F2).$  Let  {$\vp(x)$ be an eigenfunction of \eqref{eq-main-1}} and let $v(x,t)=e^{\mu t}u(x,t)$ where $u$ solves \eqref{eqn} with  nonnegative   initial data.  Then, there exists a  unique constant $\gamma^*>0$  depending on initial data such that
$$||v(x,t)-\gamma^*\vp(x)||_{C^0_x(\overline\Omega)} \to 0 \quad \mbox{as}\,\, t\to \infty.$$
\end{proposition}
\noindent {\sc Proof.}
Let us recall that   $v$ is bounded and 
$${\sup_{s\geq 1}\,||v(\cdot, \cdot+s)||_{C^{\al}_{x,t}(\overline{\Omega}\times[0,+\infty))}<+\infty\quad\mbox{for}\,\,\,\alpha>0,}$$
which  can be proved by  the Weak Harnack inequalities, \cite{W1}. Then for any sequence $\{s_n\}$, there are a subsequence $\{s_{n_k}\}$ and a function $w(x,t)$ such that
$$v(x,t+s_{n_k}) \to w(x,t)\quad\mbox{locally   in } \overline{\Omega}\times[0,+\infty)\quad \mbox{as}\,\,\,n_k\to \infty$$
and $w $ satisfies $F(D^2w)+\mu w-w_t=0$ in $\Omega\times(0,\infty).$ 
 Now  let
 $\cA $ be the set of all sequential limits of $\{v(\cdot,\cdot+s)\}_{s\geq 0}$ and let
 $$\gamma^*=\inf\{\gamma>0 : \exists\, w\in\cA \,\,\mbox{such that }\,\,w\leq \gamma\vp\,\,\,\mbox{in}\,\,\Omega\times(0,\infty)\}.$$ We  note that $0< \gamma^* <\infty$ from {\color{black}Lemma \ref{lem-asym}.} We are going to  prove that $\cA=\{\gamma^*\vp\}.$
 
First, we  show that $w\leq \gamma^*\vp$ for any $w\in\cA.$ Fix $\e>0.$   There exists  $w\in\cA$ such that
 $w\leq (\gamma^*+\e)\vp$
 by the definition of $\gamma^*.$ Then we have a sequence of functions, $\{v_n:=v(\cdot, \cdot+s_n)\},$  converging to $w$ as $s_n\to \infty,$ i.e., for a fixed $T>0,$ there is $N>0$ such that $|v_n(x,T)-w(x,T)|<\e$ for all $n>N.$ 
 Maximum principle for $e^{-\mu t}(v_n-w)$ gives us that  $|v_n(x,t)-w(x,t)|<\e$ for $\Omega\times(T,\infty).$ 
  From the Regularity Theory, 
    we have {$$||\D_x(v_n(\cdot, T+1)-w(\cdot,T+1))||_{L^{\infty}(\Omega)} \leq C\e$$} and    
hence  we deduce  $$|v_n(\cdot,T+1)-w(\cdot,T+1)|\leq C\e\vp$$ for a uniform constant $C>0$ depending on $\Omega$ and $\vp,$ i.e.,  
  $$v(x,T+1+s_n)\leq (\gamma^*+C\e)\vp(x)\quad \mbox{for large}\,\,\,s_n>0.$$
Comparison principle  implies that
$$v(x,t)=e^{\mu t}u(x,t)\leq (\gamma^*+C\e)\vp(x)  \,\,\mbox{ for} \,\,\,t\geq T+1+ s_n.$$ and also 
$$w \leq (\gamma^*+C\e)\vp \quad\mbox{for all}\,\,\,w\in\cA. $$
Since $\e$ is arbitrary and $C$ is uniform, 
$w \leq \gamma^*\vp \quad\mbox{for all}\,\,\,w\in\cA.$   

Second, we are going to show $\cA$ has only one element. Assume that $w\not\equiv \gamma^*\vp$ for some $w\in\cA.$ Then it is obvious  that $w(\cdot,0)\lneqq \gamma^*\vp$ because $u_1(x,t):=e^{-\mu t}w(x,t)$ and $u_2(x,t):=\gamma^*\vp(x) e^{-\mu t}$ solve  the same equation,
$$F(D^2u)-u_t=0\quad \mbox{in}\,\,\,\Omega\times(0,\infty).$$
Maximum principle and Hopf's Lemma imply that $u_2(x,1)-u_1(x,1)>0$ in $\Omega$ and 
$u_2(x,1)-u_1(x,1)\geq \delta\vp(x)$  for all $x\in\Omega$ for some $\delta>0,$ i.e., $w(x,1)\leq (\gamma^*-\delta)\vp(x)$ in $\Omega.$ Therefore,  we have that
$e^{\mu (t+1)}u(x,t+1)=w(x,t+1)\leq (\gamma^*-\delta)\vp(x)$ in $\Omega\times(0,\infty)$ from the comparison principle. Now, setting  $t_{n}:= s_{n}+1$  we get 
$$v(x,t+t_{n}) \to w(x,t+1)\quad\mbox{ in } \overline{\Omega}\times[0,+\infty)\quad \mbox{as}\,\,\,n\to \infty,$$
which is a contradiction to the definition of $\gamma^*.$  
Therefore we conclude that $\cA=\{\gamma^*\vp\}$ and the result follows. 
\qed

{ From the Regularity theory \label{thm-reg-1} and  the approximation lemma \ref{lem-asymp}, we obtain the following corollary.}
\begin{corollary}\label{cor-asymp} Suppose that $F$ satisfies $(F1),$ $(F2)$ and {$F$ is concave.}  Let  {$\vp(x)$ be an eigenfunction of \eqref{eq-main-1}} and let $v(x,t)=e^{\mu t}u(x,t)$ where $u$ solves \eqref{eqn} with  nonnegative   initial data.   Then we have
$$||v(x,t)-\gamma^*\vp(x)||_{C^k_x(\overline\Omega)} \to 0\quad \mbox{for some}\,\,\gamma^*>0$$
for {$k=1,2.$}
\end{corollary}


\subsection{Log-concavity}
In this subsection, we are going to study a geometric property of solutions of \eqref{eqn} and \eqref{eq-main-1} provided $\Omega$ is   convex. First, let us approximate the operator as follows.

\begin{lemma} \label{lem-F-approx}
Let $F$ satisfy   $(F1), (F2)$  and $(F3).$ 
 Then there are smooth $F_{\e}$ converging to $F$  uniformly in $Lip(\mS^{n\times n})$  satisfying $(F1), (F3)$   and 
 \begin{equation}\label{eq-F2-app}
 |DF_{\e}(z)\cdot z-F_{\e}(z)|\leq \sqrt{n}\Lambda\e.
 \end{equation}
\end{lemma}
\noindent {\sc Proof.} 
 Let  $\psi\in C_o^{\infty}(\re^{n\times n})$  be a standard mollifier with $\int\psi(z)dz=1$ and let $\psi_{\e}(z)=\frac{1}{\e^{n^2}}\psi(\frac{z}{\e^{n^2}}).$  Let us define $F_{\e}$ by  $F*\psi_\e(z).$ We note that $F_\e$ is smooth, uniformly elliptic and  concave and satisfies 
$$|F(z)-F_{\e}(z)| \leq \sqrt n\Lambda\e $$ since $F$ is uniformly elliptic. 

Now we are going  to show that for all $z,$
$$|DF_{\e}(z)\cdot z-F_{\e}(z)|\leq \sqrt n\Lambda\e.$$
Since $F$ is Lipschitz continuous with a Lipschitz constant $\sqrt n \Lambda,$  $F$ is differentiable almost everywhere from Rademacher's Theorem. Moreover, we get     $||DF||_{{\infty}}\leq \sqrt n\Lambda$ and  $$DF(z)\cdot z= F(z)\quad \mbox{a.e. z}\,\,\,\in \re^{n\times n}$$   using the fact that  
$ \displaystyle\frac{F((1+t)z)-F(z)}{t}=F(z) $ for all $z$  and  for $t>0\,$ from $ (F2).$
Then we have 
\begin{align*}
DF_{\e}(z)\cdot z-F_{\e}(z)&=\displaystyle\int\left(DF(y)\cdot z -F(y)\right)\psi_{\e}(z-y)\,dy\\
&=\displaystyle\int DF(y)\cdot\left( z - y\right)\psi_{\e}(z-y)\,dy
\end{align*}
and therefore we deduce $|DF_{\e}(z)\cdot z-F_{\e}(z)|\leq \sqrt{n}\Lambda\e.$
\qed
\begin{lemma}
\label{lem-c2f}
 Let $F$ satisfy $(F1)$, $(F2)$, and $(F3)$  and let $\Omega$ be strictly convex.  
{\color{black}Assume that $u_o\in C^0(\overline\Omega)$ be a positive initial data in $\Omega.$} If \ $\log (u_o)$ is
concave, then the solution $u(x,t)$ of \eqref{eqn}  is log-concave in
the spatial variable for all $\displaystyle0<t<\infty$, i.\,e., 
$$D^2_x\log(u(x,t))\leq 0 \quad \mbox{for }\,\,(x,t)\in\Omega\times(0,\infty).$$
\end{lemma}

\noindent {\sc Proof.}
(i) Let us    assume that  $u_o$ is smooth  in $\overline{\Omega}$ and  that $D^2\log u_o(x)\le  -c I$  in $\Omega$   for some $c>0$ 
and  approximate $F$ by  $F_{\e}$ from  Lemma \ref{lem-F-approx}.
  We also approximate   $u_o$ by $u_{\e,o}$ for small $\e>0$ such that
$${D^2 \log u_{\e,o}\leq0\,\,\,\mbox{in}\,\,\,\Omega,\,\,\,\,\,F_{\e}(D^2u_{\e,o})=0\quad\mbox{on}\quad\p\Omega.}$$

Then there is the positive smooth solution $u_{\e}$ of \eqref{eqn}
with an operator { $F_{\e}(\cdot)-F_{\e}(0)$ } and an initial data $u_{\e,o}$. Let us  put $g(x,t)=\log u_{\e}(x,t)$, which is
finite and smooth for $x\in \Omega$ and takes the value $g=-\infty$
on  $\partial\Omega\times (0,\infty)$. It also
satisfies the equation
\begin{equation*}
\pa_t g =e^{-g}F_{\e}\left(e^{g}\left(D^2g+\D g \D g^T\right)\right)-{e^{-g}F_{\e}(0)}\quad\mbox{in}\,\,\,\Omega\times(0,+\infty).
\end{equation*}

First, let us consider  a domain $\Omega\times(0,T)$ for  $T>0.$ 
To estimate the maximum of its second derivatives, for small $\delta>0,$
 consider the function $Z$ defined as 
$$
Z(t)=\sup_{y\in \Omega}\sup_{|e_{\beta}|=1} \, g_{\beta\beta}(y,t)+\psi(t)
,$$
where  $e_{\beta}\in S^{n-1}$ and 
 $\psi(t):=-\delta\tan(2K\sqrt{\delta}\,t).$  The constant  $K>0$ independent of $\e>0$ and $\delta>0$  will be chosen later. 
Now, 
 let us assume there exists $ t_o\in\left[0,\min\left(\frac{\pi}{4K\sqrt{\delta}},T\right)\right]$ such that $$Z(t)=\sup_{y\in \Omega}\sup_{|e_{\beta}|=1} \, g_{\beta\beta}(y,t)+\psi(t)=0\quad \mbox{at}\,\,t=t_o.$$ We may assume that 
 $$Z(t_o)=g_{\al\al}(x_o,t_o)+\psi(t_o)=0$$ for some direction $e_\al$ and $x_o\in\overline{\Omega}.$
Then $e_\alpha$ is an
eigen-direction of the symmetric matrix $D^2g(x_o,t_o)$ which means
that, using orthonormal coordinates in which $e_\alpha$ is taken as
one of the coordinate axes,  $ g_{,\alpha\beta}$ is zero at
$(x_o,t_o)$ for $\beta\neq\alpha.$  We note that $Z(0)<0$ from the assumption.

 Then, we claim that
$$
g_{,\alpha\alpha}(x,t_o)=\frac{u_{\e}\,u_{{\e},\alpha\alpha}- u_{{\e},\alpha}^2}{u_{\e}^2}\ra
-\infty\quad\mbox{as}\,\,\,x\in \Omega \ra \partial\Omega.
$$
This holds when  $e_\al$ is not a tangential direction,
since $ \partial \Omega$ is
smooth, $|D^2 u_{\e}|$ is bounded   and  $|\D u_{\e}|>0$ on $\partial\Omega$ by Hopf's lemma. 
For a tangential direction $e_\al,$
we take a coordinate system such that $x_o=0$
and that the tangent plane is $x_n=0.$ Let the boundary be given
locally by the equation $x_n=f(x')$, and  $x'=(x_1,\cdots,x_{n-1}).$ We introduce the change of variables
$$
y_i=x_i \,\,(i=1,\cdots, n-1), \quad y_n=x_n-f(x'), \quad v(y,t)=u_{\e}(x,t) .
$$
Then along tangent directions $e_\al$ we have
$$
u_{{\e},\alpha\alpha}(x,t)=v_{,\alpha\alpha}(y,t)-2v_{,n\alpha}(y,t)\,f_{,\alpha}(x')+v_{,nn}(y,t)\,(f_{,\alpha}(x'))^2
-v_{,n}(y)\,f_{,\alpha\alpha}(x').
$$
Using the fact that $\p_{jj}v(0,t)=0$  from the boundary condition   and  $f_{j}(0)=0$ for  $j=1,\cdots, n-1,$ we obtain
$$
u_{{\e},\alpha\alpha}(0,t_o)= -v_{n}(0)\,f_{\alpha\alpha}(0)<0,
$$
for a tangential vector $e_\al.$
We note that $f_{,\alpha\alpha}(0)>0$ since  $\Omega$ is strictly convex. 
Thus   $g_{,\alpha\alpha}(x,t_o)$ tends to $-\infty$ as $x\in\Omega$ goes to $\partial  \Omega$.
And from the uniform global $C^{2,\be}$ estimate of $u_{\e}$,  there is a small $\eta>0$ independent of $\e,\delta$ such that   {$g_{\alpha\alpha}(x,t)<-10$ for $x\in \Omega\backslash\Omega_{(-\eta)}\times(0,T),$}
where $\Omega_{(-\eta)}=\{x\in\Omega: d(x,\partial\Omega)>\eta\}$.
{\color{black} So we deduce that the maximum of $Z$ can only  be achieved at an
interior point  $x_o\in\Omega_{(-\eta)}$.}

Next, we   look at the evolution equation of
$g_{\alpha\alpha}(x,t),$ which is given by the equation as below
\begin{align*}
g_{\al\al,t}&=F_{ij}\cdot(D_{ij}g_{\al\al}+D_ig_{\al\al}D_jg +D_igD_jg_{\al\al}+ 2D_ig_\al D_jg_\al  ) \\
&\,\,\,+(g_{\al}^2 -g_{\al\al})\left\{ e^{-g}F\left(e^g\left(D_{ij}g+D_igD_jg\right)\right)- F_{ij}\cdot\left(D_{ij}g+D_igD_jg\right)\right\}\\
&\,\,\,+e^{-g}F_{ij,kl}\cdot\left(e^{g}\left(D_{ij}g+D_ig D_j g\right)\right)_{\al}\left(e^{g}\left(D_{kl}g+D_{k}g D_l g\right)\right)_{\al}\\
&\,\,\,{-(g_{\al}^2 -g_{\al\al}) e^{-g}F_{\e}(0)}
\end{align*}
where $F_{ij}={F_{\e,}}_{ij}\left(e^{g}\left(D_{ij}g+D_igD_jg\right)\right).$
Since $F_\e$ satisfies $(F1),$ concavity and \eqref{eq-F2-app}, it follows that
\begin{align*}
g_{\al\al,t} \,
\leq F_{ij}\cdot(D_{ij}g_{\al\al}+D_ig_{\al\al}D_jg +D_igD_jg_{\al\al}+ 2D_ig_\al D_jg_\al  ) +2\sqrt{n}\Lambda\e e^{-g}|g_{\al}^2 -g_{\al\al}|.
\end{align*}
At the point of maximum $(0,t_o),$ we see that $$g_{\al\al}=-\psi \geq 0, \,\,\D_x g_{\,\alpha\alpha}=0 , \,\,D_x^2
g_{\,\alpha\alpha}\le 0$$ as well as $g_{\,\alpha\beta}=0$ for
$\beta\ne \alpha.$
Thus   we  get  at the point of maximum $(0,t_o),$ 
\begin{align*}
g_{\al\al,t}&\leq F_{ij}\cdot(D_{ij}g_{\al\al}+D_ig_{\al\al}D_jg +D_igD_jg_{\al\al}+ 2D_ig_\al D_jg_\al  ) +2\sqrt n\Lambda\e e^{-g}|g_{\al}^2 -g_{\al\al}|\\
&\leq 2F_{\al\al} g_{,\alpha\alpha}^2  +2\sqrt n\Lambda\e e^{-g}(g_{\al}^2 +g_{\al\al})\\
&\leq2\Lambda g_{,\alpha\alpha}^2 +\e2\sqrt n\Lambda \frac{|u_{\e,\alpha\alpha}|}{u_{\e}}.
\end{align*}

On the other hand, when the supremum of  $Z(t)-\psi(t)=\sup_{y\in \Omega}\sup_{|e_{\beta}|=1} \, g_{\,\beta\beta}(y,t)$ is achieved at   a point $x(t)\in\Omega$  with a unit vector $e_{\beta(t)}$ at each time t,  we  check  that $g_{\,\be(t)\,\beta'(t)}=0$  and $\D_x g_{\,\be(t)\,\be(t)}=0$ at the point $(x(t),t).$
Therefore, we have  at the maximum point  $(0,t_o),$ 
\begin{align*}
0\leq Z'(t_o)&=g_{\al\al,t}+\psi_t \\&\leq \psi_t+2\Lambda\psi^2+\e2\sqrt n\Lambda\frac{|u_{\e,\alpha\alpha}|}{u_\e}\leq \psi_t+K(\psi^2 +\e) ,
\end{align*}
when we select a uniform number $K>0$ bigger that {  $C(\Lambda,n)\left(1+\displaystyle\max_{\Omega_{(-\eta)}\times(0,T)} \frac{|D^2u_{\e}|}{u_{\e}}\right).$  } 
Now, it is easy to  check that $$ \psi_t+K( \psi^2+\e)<\frac{2K(-\delta^{3/2}+\delta^2)}{\cos(2K \sqrt\delta t)}<0$$ for $0<\e<<\delta$ and for $2K \sqrt\delta t<\frac{\pi}{2},$ which
implies  a contradiction. 
 Therefore, we obtain  $$\sup_{y\in \Omega}\sup_{|e_{\al}|=1}  \partial_{\alpha\alpha}\log (u_{\e})(y,t)< -\psi(t)=\delta\tan(2K\sqrt\delta t)\leq \delta$$
for $\displaystyle0<t<\min\left(\frac{\pi}{8K\sqrt{\delta}},T\right)$ and for $0<\e<<\delta$ from  the uniform interior { $C^{2,\be}$-estimates of $u_{\e}$ in $\Omega_{(-\eta)}\times(0,T).$}   
Letting $\delta\ra 0$  we conclude that
 $$\partial_{\alpha\alpha}\log (u)\leq 0 \quad\mbox{in}\,\,\,\,\Omega\times(0,T).$$
Therefore $u(x,t)$ is log-concave with respect to $x$ in $\Omega\times(0,\infty)$ since $T$ is arbitrary.

(ii) The proof in the general case uses a density argument which
is more or less standard. Briefly, if $u_o$ is not smooth and
strictly log-concave, we  perform a mollification to obtain
an approximating  sequence $u_{oj}$ of smooth and log-concave
functions. 
To make $u_{oj}$ strictly
log-concave we may  put for instance,
$$
{\tilde u}_{oj}(x)= u_{oj}(x)\,\mbox{exp\,}(-c_j |x|^2)
$$
for some $c_j >0$, $c_j \to 0$ as $j\to\infty.$ From (i), we get the
result for $\tilde u_j$, the solution of the problem with data
${\tilde u}_{oj}.$
Uniform H\"{o}lder regularity 
let us take a subsequence $\tilde u_j$ converging uniformly to
$u$ in each compact subset and then uniform convergence on each compact subset will preserve the sign of the second difference  in the limit.
\qed

\begin{corollary}\label{cor-conc}
Let  $F$ satisfy $(F1)$,$(F2)$ and $(F3)$ and let { $\Omega$ be  convex.} 
If  $u_o$ is log-concave, so is the viscosity solution $u(x,t).$
\end{corollary}

\begin{remark}\label{rmk-conc}
{\rm We note that any concave function in a convex
domain $\Omega$ is  log-concave. On the other hand, it is well-known
that the distance function $ \mbox{dist}(x,\partial \Omega)$ is
concave for a  convex domain, so the lemma is not void.}
\end{remark}
\begin{remark}\label{rmk2-conc}
\item

\begin{enumerate}
\item Let $\sigma_k(D^2 u)=\sum_{i_1<\cdots<i_k}\lambda_{i_1}\cdots\lambda_{i_k}$ for the eigenvalues $\lambda_1\leq \cdots \leq\lambda_n$ of $D^2 u$. $F(D^2 u)=\sigma_k(D^2 u)^{\frac{1}{k}} $ satisfies  the conditions  $(F2)$ and $(F3)$. 
\item If  a differentiable operator $F$ satisfies (F2),  then $F$ is  linear. If $F$ is also uniformly elliptic, then $F$ becomes Laplacian after suitable trasformation. 
\end{enumerate}

\end{remark}
\begin{corollary}[Log-concavity]\label{cor-log-con}
Let  $F$ satisfy $(F1)$,$(F2)$ and $(F3)$ and let
 $\Omega$ be convex. Then, the stationary profile $\vp(x)$   is
log-concave, i.\,e., \ $D^2\log(\vp(x))\leq 0$.
\end{corollary}

\noindent {\sc Proof.}
Take  the distance function as an initial data  of parabolic flow,\eqref{eqn}. 
 Then Corollary \ref{cor-conc} yields that for $x,y\in\Omega,$ 
$$2\left(\log u(x,t)+\log u(y,t)\right)-\log u\left(\frac{x+y}{2},t\right) \leq 0.$$
From the
asymptotic result, Proposition \ref{lem-asymp}, we have the uniform
convergence
$$
||e^{\mu
t}u(x,t)-\gamma^*\vp(x)||_{{C^0_x}(\overline\Omega)}\ra 0\quad\mbox{as}\quad t\to\infty
$$ and hence the result follows.
\qed

For a differentiable operator,  the foregoing is a classical result, \cite{LV2}  for a  domain which is smooth and strictly convex.

\begin{lemma}[Strict log-concavity]\label{lem-sc-1}
\label{lem3}
Suppose that $F$ satisfies (F1), (F2) and (F3) and is differentiable and that  $\Omega$ is smooth and strictly convex. Then the positive eigenfunction $\vp$ of \eqref{eq-main-1} is strictly log concave, i.e., there exists a constant $c_1>0$ such that
\begin{equation*}
D^2(\log\,
\vp)\leq -c_1\,{\bf I}.
\end{equation*}
\end{lemma}



\begin{theorem}[Eventual log-concavity]
\label{theor1}  We assume the same hypothesis as Lemma \ref{lem-sc-1}. Let $u_o$ be a nonnegative initial
function.  
Then, the solution $u(x,t)$ of   {\rm \eqref{eqn}}
is strictly log-concave in the spatial variable for large $t>0,$
i.e., for every $\varepsilon>0$ there is
$t_o=t_o(u_o,\varepsilon)$ such that
\begin{equation*}
 D^2\log(u(x,t))\leq
-(c_1-\varepsilon)\,{\bf I} \qquad \mbox{for all } t\ge t_0,
\end{equation*}
where $c_1>0$ is the constant of Lemma {\rm \ref{lem3}}.
\end{theorem}



\section{Degenerate Parabolic Fully Nonlinear Equation}\label{sect-pme}

In this section, we consider the solution $u(x,t)$ of the fully nonlinear degenerate parabolic equation
\begin{equation}\label{eqn-pme}
\left\{
\begin{array}{ll}
u_t(x,t)=F(D^2 u^m(x,t)) \qquad &\text{in $Q_T=\Omega\times (0,T)$},\,\,\,m > 1,\\
u(x,0)=u_o(x), & \\
u(x,t)=0\quad &\text{on $x\in \partial\Omega$},
\end{array}\right. \qquad \qquad
\end{equation}
where $\Omega$ is a bounded domain of $ \re^n$ with a smooth
boundary. We assume that $u_o^m\,$ belongs to {$$\cC(\overline\Omega):=\{ h\in C^o(\overline\Omega) \,|\, c_o\,\mbox{dist}(x,\p\Omega)\,\leq h(x) \leq C_o\,\mbox{dist}(x,\p\Omega)\,\,\mbox{for some } 0<c_o\leq C_o<+\infty\}.$$}

We  define  $w:=u^m,$ then $w$ satisfies 
 \begin{equation}\label{eqn-pme-m}
 \left\{
\begin{array}{ll}
m w^{1-\frac{1}{m}}F(D^2w)-w_t=0 \quad &\text{in $Q_T=\Omega\times (0,T)$},\,\,\,m > 1,\\
w(x,0)=w_o(x)=u_o^m(x)\in \cC(\overline\Omega), & \\
w(x,t)=0\quad &\text{on $x\in \partial\Omega$}.
\end{array}\right. \qquad \qquad
\end{equation}
We also introduce the pressure in the form $v=\frac{m}{m-1} u^{m-1}.$ If $F$ satisfies $(F2),$ the pressure  $v$ solves 
\begin{equation}\label{eq-pres}
\left\{
\begin{array}{ll}
v_t = F((m-1)vD^2v + DvDv^T) &\text{in $Q_T=\Omega\times (0,T)$},\,\,\,m > 1,\\
v(x,0)=v_o(x)=\frac{m}{m-1}u_o^{m-1},\,\,&\,\,   \\
v(x,t)=0\quad &\text{on $x\in \partial\Omega$}.
\end{array}\right. \qquad \qquad
\end{equation}

Before studying   asymptotic behaviors of degenerate parabolic flows, let us state  the regularity of the solution.

\begin{proposition} [Regularity for $m>1$]\label{prop-reg}
Let $F$ satisfy $(F1),(F2)$ 
and let $u$ be the solution of $\eqref{eqn-pme}.$
 \begin{enumerate}
\item If  $u_o$ is nonzero and  nonnegative, then 
\begin{enumerate}[(i)]
\item 
there is a time $t_o=t_o(u_o,\Omega) > 0\,\,$ such that
$$u(x,t)  >0\qquad\mbox{in}\,\,\,\Omega\times(t_o,\infty)$$ for a uniform constant $t_o=t_o(\lambda,\Lambda,u_o)>0.$\\

\item   $0 \leq u(x,t) \leq {C_o}{t^{-\frac{1}{m-1}}} \,\mbox{dist}\,(x,\p \Omega)^{\frac{1}{m}}\,\,\,\mbox{in}\,\,\, \Omega\times(0,\infty).$
\end{enumerate}

\item If $u_o$ is an initial data in $\cC(\overline\Omega),$ then  
\begin{enumerate}[(i)]
\item 
 we have 
$${c_o}{(t+\tau_1)^{-\frac{1}{m-1}}} \,\mbox{dist}\,(x,\p \Omega)^{\frac{1}{m}} \leq u(x,t) \leq {C_o}{(t+\tau_2)^{-\frac{1}{m-1}}} \,\mbox{dist}\,(x,\p \Omega)^{\frac{1}{m}}\,\,\,\mbox{ in}\,\,\, \Omega\times(0,\infty)$$
for some constant $\tau_1,\tau_2$ depending on $u_o.$  Moreover, for  $Q_T=\Omega\times{[s, T]}, (0<s<T),$  
\begin{enumerate}[(a)]
\item $u$ is of $C^{1,\beta}(Q_T)$ for some $0<\beta<1,$
\item   $u$ is of $C^{1,1}(Q_T)$  if $F$ is concave or convex, 
\item   $u$ is of $C^{2,\beta}(Q_T)$  for some $0<\beta<1$ if $F$ is concave or convex and $u\in C^{1,1}$, 
\item $u$ is of $C^{\infty}(Q_T)$   if $F$ is $C^{\infty}$ and $u\in C^{2,\beta}$.  \\
\end{enumerate}

\item  $u$ is of $C_x^{0,\frac{1}{m}}(\overline{\Omega}\times [s,T])\cap C^{1,\beta}(\Omega\times [s,T]) $ for some $0<\beta<1.$ 
\end{enumerate}
\end{enumerate}
 \end{proposition}
 
\noindent {\sc Proof.}
(1)  For $c > 0$, let
$$V(x,t) = t^{-\al}{\left(c-k\frac{|x|^2}{t^{\be}}\right)}_+ \,\,, $$
$\mbox{where}\,\,\, \al=\frac{n(m-1)\Lambda}{2\lambda + n(m-1)\Lambda},\, \,
\be=\frac{2\lambda}{2\lambda + n(m-1)\Lambda},\,\,\mbox{and}\,\,\,\,k=\frac{1}{2(2\lambda + n(m-1)\Lambda)} .$ Then we can check 
\begin{align*}
&F((m-1)VD^2V + DV DV^T) - V_t\, \\&\geq\, \cM^-((m-1)VD^2V) + \cM^-(DV DV^T) - V_t
=0\quad\mbox{in}\,\,\,\{V>0\}
\end{align*}
and hence $V$ is a sub-solution of \eqref{eq-pres} as long as $\supp(V)\subset\overline\Omega$.

We define $\,\overline{U}(x,t) = \left(\frac{m-1}{m}V(x,t)\right)^{\frac{1}{m-1}}
=\left(\frac{m-1}{m}\right)^{\frac{1}{m-1}}
t^{-\al/(m-1)}\left(c-k\frac{|x|^2}{t^{\be}}\right)_+^{\frac{1}{m-1}},\,\,  $
and hence  $\overline{U}\,$ is a subsolution of $\eqref{eqn-pme}$ in $\supp(\overline{U})$ as long as $\supp(\overline{U})\subset\overline\Omega$.
We note that the support of $\overline{U}$ is compact and expands in time. So the  previous argument in Lemma \ref{lem-asym} gives the result that $u$ is positive for large time $t.$

(ii)  To get the upper bound, we are going to show that
$$u(x,t)\leq f(x) t^{-\frac{1}{m-1}}\,\,\, \,\,\mbox{in}\,\,\, \Omega\times(0,\infty),$$
where $f$ is the solution of $\eqref{eq-main-1-pme}.$ 
Define  $u_{o,\e}:=(u_o-\e)_+=\max(u_o-\e, 0)$ for $\e>0$ and  let $u_{\e}$ be the solution of \eqref{eqn-pme} with initial data $u_{o,\e}.$ We choose $\tau_\e>0$ converging to $0$ as $\e\to0$ such that $u_{o,\e}(x)\leq f(x) (\tau_\e)^{-\frac{1}{m-1}}.$ Comparison principle yields that  $$u_{\e}(x,t)\leq f(x) (\tau_\e+t)^{-\frac{1}{m-1}}\leq f(x)t^{-\frac{1}{m-1}}$$ in $\Omega\times(0,\infty)$ since $ f(x) (\tau+t)^{-\frac{1}{m-1}}$  is a similarity solution for any $\tau>0.$ 
 
From the comparison principle, $u_{\e}$ is nondecreasing  as $\e$ decreases and 
$$u_{\e_o}\leq u_{\e}\leq u\leq M:=\max_{\Omega} u_o$$ if $\e<\e_o$ for  any $\e_o>0.$ 
Then for  each compact subset $K$ of $\Omega\times(0,\infty),$ $w_\e:=u^m_\e $ satisfies   a uniformly parabolic equation, $w^{1-\frac{1}{m}}F(D^2w)-w_t=0,$ and uniform parabolic estimates tell us that $w_{\e}\to \tilde w$ as $\e\to 0$ in $K$ for some locally H\"older continuous function $\tilde w,$ which  is the solution of $\eqref{eqn-pme-m}.$ Therefore, we obtain
$$u(x,t)\leq f(x) t^{-\frac{1}{m-1}} \,\,\mbox{in}\,\, \Omega\times(0,\infty)$$ and hence $0 \leq u(x,t) \leq {C_o}{t^{-\frac{1}{m-1}}} \,\,\mbox{dist}\,(x,\p \Omega)^{\frac{1}{m}}$ in $\Omega\times(0,\infty)$
since    $\displaystyle\inf_{\p\Omega}|\D_x f^m| >0.$

(2) (i) We choose $\tau_1>0, \tau_2>0\,$ such that 
$$f\cdot\tau_1^{-\frac{1}{m-1}} \leq u(\cdot, 0) \leq f\cdot\tau_2^{-\frac{1}{m-1}}$$
because $u_o^m\in \cC(\overline\Omega). $
Since $f(x)(\tau_i+t)^{-\frac{1}{m-1}}$ is a solution of $\eqref{eqn-pme},$  the comparison principle implies 
$$f\cdot(\tau_1+t)^{-\frac{1}{m-1}} \leq u(\cdot,t) \leq f\cdot(\tau_2+t)^{-\frac{1}{m-1}}.$$
Thus the first  result  comes from the gradient estimate of the positive  eigenfunction on the boundary.  
On the other hand, for each compact subsets $K\Subset\tilde K\Subset \Omega,$  there exist   $0<c_o\leq C_o <+\infty$ such that
 $$0<c_o \leq w(x,s)=u^m(x,s)\leq\ C_o <+\infty\,\,\,\mbox{in}\,\,\, \tilde K\times[s/2,T],$$
which  means that the operator $w^{1-\frac{1}{m}}F(\cdot)$ becomes uniformly elliptic in $ \tilde K\times[s/2,T].$ So the estimates  follow from Theorem \ref{thm-reg-1}.

(ii) We  use  the fact (i) and scaling property to prove the H\"older regularity on the boundary.
In fact, since 
 we have a linear growth of $w=u^m$ away from the boundary: let $\delta_o>0$ be a constant such that $B_{\delta_o}(x)\subset\Omega$ for $\mbox{dist}\,\,(x,\p\Omega)>\delta_o.$ 
 For  $x_o\in\Omega$ such that $\mbox{dist}\,\,(x_o,\p\Omega)<\delta_o,$ we  set $\,\mbox{dist}\,(x_o,\p\Omega)=2\sigma.$ According to (i), it follows that
\begin{equation}\label{eqn-pme-bd-sc}
c_o\sigma < w(x_o,t)=u^m(x_o,t) < C_o \sigma,\,\,\,\mbox{for}\,\,\,t\in[s/2,T],
\end{equation}
where $|x-x_o|=\,\mbox{dist}\,(x,\p\Omega) =2\sigma<\delta_o.$ 

Now we scale $w$ linearly with the distance $\sigma$ to the boundary so the scaled function $\tilde w$  has a value of order one. Then  $\tilde{w}$ will satisfy a uniformly parabolic equation and have a uniform gradient estimate.
 Define $\tilde{w}(\tilde x,\tilde t)= w^{\sigma}(\tilde{x},\tilde{t}) := \frac{1}{\sigma}w(x_o+\sigma \tilde{x},\sigma^{1+1/m}\tilde{t}).$ From scaling property, $\tilde{w} $ satisfies $\tilde{w}^{1-\frac{1}{m}}\tilde F(D^2\tilde w)-\tilde{w}_t=0$ 
 for an elliptic operator  $\tilde F(\cdot)=\sigma F\left(\frac{\cdot}{\sigma}\right)$ with the same ellipticity constants $\lambda,\Lambda$ and this transform sends
$\{x\in\Omega : \,||x-x_o||\,=\sigma\}$ to $\{\tilde x : \, ||\tilde x||=1\}$.
Thus $\eqref{eqn-pme-bd-sc}$ implies that
$$c_o< \tilde{w}(\tilde{x},\tilde{t})< C_o \,\,\,\,\mbox{for}\,\,(\tilde x,\tilde t)\in B_1(0)\times[\sigma^{-1-1/m} s/2,\sigma^{-1-1/m}T],$$
and then we have $$|\D w(x_o,t)|=|\D \tilde{w}(0,\tilde{t})|< C\,\,\,\mbox{for}\,\,\,t\in[s,T]$$
 from uniform gradient estimates for uniformly parabolic equations.( We refer to \cite{L},\cite{W1}.)

On the region $K:= \{x\in\Omega : \mbox{dist}\,({x},\p{\Omega})\,\geq\frac{\delta_o}{2}\},$  $u$ is positive so we have $u\geq c_o\,\,\mbox{in}\,\, K\times[s/2,T]$ for some constant $c_o>0.$
Then  the operator is uniformly parabolic in $K\times[s/2,T]$ and hence  we also have
$${|\D w (x,t)| < C||w||_{L^{\infty}(K\times[s/2,T])}}$$
for  $\mbox{dist}\,({x},\p{\Omega})\,\geq\frac{3}{4}\delta_0$ and ${t \in[s,T]}$ from the regularity theory of the uniformly parabolic equations. Therefore, $w=u^m$ is of $C^{0,1}(\overline{\Omega}\times [s,T]).$ 
\qed

\subsection{Asymptotic Behavior}
First, we are going to show Aronson-B\'enilan  inequality for the degenerate fully nonlinear equation with $m>1$,
which tells us almost monotonicity of parabolic flows as $t\ra \infty$.

\begin{lemma}[Aronson-Benilan inequality]\label{lem-ab}
{Suppose that $F$ satisfies $(F1),$ $(F3)$ and $F(0)=0.$}
Let $u$  be the solution of \eqref{eqn-pme}  with  initial data $u_o^m\in\cC(\overline\Omega)$ and let $v=u^{m-1}.$ Then
we have for large $C=C(m)>0,$
\begin{equation}\label{ineq-ab}
u_t\geq -C\frac{u}{t}\qquad \mbox{and}\qquad v_t\geq -C\frac{v}{t}\qquad\mbox{for}\,\,\,t>0.
\end{equation}
\end{lemma}
\noindent {\sc Proof.}
(i) First, let us also assume that $F$ is of $C^1.$   Let $w:=u^m$ and $w$ solves \eqref{eqn-pme-m}.  Let $\delta>0$ and $\e>0$ and let $C$ be a positive constant bigger than $\frac{m}{m-1}.$   We 
can select $-\delta<\tau_{\e,\delta}<0$ so that $\displaystyle w_t+{C}\frac{w+\e}{t+\tau_{\e,\delta}} >0$ at $t=\delta$ because $w_t=w=0 $ on $\p\Omega\times(0,\infty).$

Define $Z(t):=\displaystyle\inf_{x\in\Omega}\left(w_t+C\frac{w+{\e}}{t+\tau_{\e,\delta}}\right).$  We note that  $Z(\delta) > 0.$  
 From the concavity of $F$ and $F(0)=0,$ the function  $w$ satisfies 
$$mw^{1-1/m}F_{ij}(D^2w)D_{ij}w-w_t\leq 0.$$
Let $t_o\in(\delta,\infty)$  be the first time such that $Z(t_o)=0$ and then we have that  $w_t(x_o,t_o)<0,$ 
and that $\displaystyle{w_t^2}=C^2\displaystyle\frac{(w+{\e})^2}{(t+\tau_{\e,\delta})^2}>0$ at the minimum point $(x_o,t_o)\in\Omega\times\{t_o\}.$  Indeed,  the minimum point  $x_o$ is   interior in $\Omega$ because $\displaystyle\left(w_t+C\frac{w+{\e}}{t+\tau_{\e,\delta}}\right)>0$ on $\p\Omega.$  At the minimum point, we have
\begin{equation*}
\begin{split}
Z_t&=\left(w_t+C\frac{w+{\e}}{t+\tau_{\e,\delta}}\right)_t\\
&=\left(1-\frac{1}{m}\right)mw^{-1/m}F(D^2w)w_t+mw^{1-1/m}F_{ij}(D^2w)D_{ij}w_t+C\frac{w_t}{t+\tau_{\e,\delta}}
-C\frac{w+{\e}}{(t+\tau_{\e,\delta})^2}\\
&=\left(1-\frac{1}{m}\right)\frac{(w_t)^2}{w}+mw^{1-1/m}F_{ij}D_{ij}\left(w_t+C\frac{w+{\e}}{t+\tau_{\e,\delta}}
\right)\\
&\,\,\,\,-mw^{1-1/m}F_{ij}D_{ij}\left(C\frac{w+{\e}}{t+\tau_{\e,\delta}}
\right)+C\frac{w_t}{t+\tau_{\e,\delta}}
-C\frac{w+{\e}}{(t+\tau_{\e,\delta})^2}\\
&\geq\left(1-\frac{1}{m}\right)\frac{(w_t)^2}{w}-C\frac{w_t}{t+\tau_{\e,\delta}}+C\frac{w_t}{t+\tau_{\e,\delta}}
-C\frac{w+{\e}}{(t+\tau_{\e,\delta})^2}\\
&=\left(1-\frac{1}{m}\right)C^2\displaystyle\frac{(w+\e)^2}{w(t+\tau_{\e,\delta})^2}-C\frac{w+\e}{(t+\tau_{\e,\delta})^2}\geq C\displaystyle\frac{w+\e}{(t+\tau_{\e,\delta})^2}\left(\frac{m-1}{m}C-1\right)>0,
\end{split}
\end{equation*}
which is a contradiction. 
Therefore we have
$w_t > -C\displaystyle\frac{w+\e}{t+\tau_{\e,\delta}}\geq -C\displaystyle\frac{w+\e}{t-\delta}$ for $t>\delta.$  
 Since  $\e,\delta>0$ are arbitrary, we deduce that $tw_t+Cw\geq 0$    for $\Omega\times(0,\infty)$ and hence $u_t\geq -C\displaystyle\frac{u}{t}$ for $t>0.$  

(ii) In general, let us approximate $F(\cdot)$ by smooth $F_{\e}(\cdot).$ 
 Let ${u^{\e}}$ be the solution of \eqref{eqn-pme} with the operator $F_{\e}$ and with the same initial data and let  $u^{\pm}$  be the solution of \eqref{eqn-pme} with the operator  $\cM^{\pm}.$  Let us define $w^{\e}:=({u^{\e}})^m.$ From Comparison principle,  it follows that $$0<u^{-} \leq u^{\e}\leq u^{+}\leq\|u_o\|_{L^\infty(\Omega)}\quad \mbox{in} \,\,\,\Omega\times(0,\infty),$$ which implies that   $w^{\e}$ solves the  uniformly parabolic equation  in each compact subset of $\Omega\times(0,\infty).$
 Then, $w^{\e}$ and $w_t^{\e}$ converge uniformly to  $w$ and $w_t$, respectively, in each compact subset  of $\Omega\times(0,\infty)$  from the regularity theory. 
 Therefore we conclude  that  $w_t\geq -C\displaystyle\frac{w}{t}$ 
for large $C=C(m)>0$ and hence $\eqref{ineq-ab}$ holds by direct calculations.
 \qed
 


\begin{proposition}[Approximation]\label{prop-pme-app} Suppose that $F$ satisfies $(F1)$ and $(F2)$.
Let $u$ be the solution of $\eqref{eqn-pme}$ with initial data $u_o^m\in\cC(\overline\Omega).$  Set $U(x,t) := \frac{f(x)}{(1+t)^{\frac{1}{m-1}}},\,\,$ where $f\,\,$ solves\begin{equation}
\begin{cases}
 -F(D^2 \ f^m(x)) = \frac{1}{m-1}f(x)\quad\text{in $\Omega$},\,\, m>1,\\
 f = 0\quad\text{on $\partial\Omega,$}\\
 f > 0\quad\text{in $\Omega$}.
\end{cases}
\label{eq-main-1-pme-1}
\end{equation} Then,  we have
$$\lim_{t\to\infty} t^{\frac{1}{m-1}}|u(x,t) - U(x,t)| \to 0\,\,\,\,\mbox{uniformly in }\,\,\, \Omega,$$
 and  there exists  $ t_o>0\,\,$ such that $u^m\,\,$ is $C^1\,\,$ up to the boundary
and $0 < c_o < t^{\frac{m}{m-1}}|\D u^m(x,t)| < C_o \,\,$ for $x\in\p\Omega,$ where  $c_o\,$ and $C_o$  depend  on $u_o$ and $\Omega.$ 
\end{proposition}
\noindent {\sc Proof.}
(i) In the proof of (2) at Proposition \ref{prop-reg}, we have
$$f\cdot(\tau_1+t)^{-\frac{1}{m-1}} \leq u(\cdot,t) \leq f\cdot(\tau_2+t)^{-\frac{1}{m-1}}$$
since $u_o^m\in\cC(\overline\Omega).$ 
Then, we obtain
$$  t^{\frac{1}{m-1}}|u - U| \leq f\cdot\left(  \frac{t^{\frac{1}{m-1}}}{(\tau_2+t)^{\frac{1}{m-1}}}- \frac{t^{\frac{1}{m-1}}}{(\tau_1+t)^{\frac{1}{m-1}}}\right)\to 0\quad\mbox{
uniformly}\,\,\,\,\mbox{as}\,\,t\to\infty.$$

 (ii) From (i), $w=u^m$ satisfies 
$$\phi\cdot(\tau_1+t)^{-\frac{m}{m-1}} \leq w \leq \phi\cdot(\tau_2+t)^{-\frac{m}{m-1}} \quad\mbox{in $\Omega\times[0,\infty),$}$$ where $\phi=f^m$ is the solution of \eqref{NLEV}.
 From Hopf's Lemma for $\phi,$  we have $$  \frac{c_1}{(1+t)^{\frac{m}{m-1}}}\,\mbox{dist}\,(x,\p\Omega)\,\,\leq w(x,t) \leq \,\frac{c_2}{(1+t)^{\frac{m}{m-1}}}\,\mbox{dist}\,(x,\p\Omega)\quad\mbox{in}\,\,\Omega\times{[0,\infty)}$$ and 
 $$  \frac{c_1}{t^{\frac{m}{m-1}}}\,\mbox{dist}\,(x,\p\Omega)\,\,\leq w(x,t) \leq \,\frac{c_2}{t^{\frac{m}{m-1}}}\,\mbox{dist}\,(x,\p\Omega)\quad\mbox{in}\,\,\Omega\times{[1,\infty)}.$$
We follow a similar argument as (2),(ii) at Proposition \ref{prop-reg} and use scaling property for porous medium equation to estimate 
$$  |\D w(x,t)| \leq \,\frac{C_o}{t^{\frac{m}{m-1}}}\,\,\,\,\mbox{for}\,\,(x,t)\in\overline\Omega\times{[1,\infty)}$$
 for  some $0<C_o<+\infty.$
 Moreover,  we have that
$$\frac{c_o}{t^{\frac{m}{m-1}}} \leq  |\D w(x,t)| \leq \,\frac{C_o}{t^{\frac{m}{m-1}}}\,\,\,\,\mbox{for}\,\,(x,t)\in\p\Omega\times[1,\infty)$$
 for some $0<c_o\leq C_o<+\infty,$
which means that $0 < c_o < t^{\frac{m}{m-1}}|\D u^m(x,t)| < C_o \,\,$ for $x\in\p\Omega$ and  for large  $t>t_o.$ 
 \qed

\begin{remark}{\rm
If we set $z(x,t):=\displaystyle t^{\frac{1}{m-1}}u(x,t)$, the renormalized function, the estimate in Lemma \ref{lem-ab}
 holds to $z.$ In fact, we have
\begin{align*}
z_t=\displaystyle t^{\frac{1}{m-1}}\left(u_t+\frac{1}{m-1}\frac{u}{t}\right)\geq \left(-C+\frac{1}{m-1}\right)\frac{z}{t}.
\end{align*}}
\end{remark}

\begin{corollary}\label{cor-pme-app}
\item Under the same condition of Proposition \ref{prop-pme-app},
          \begin{equation}\label{eq-rn-d}
           z(x,t)=t^{\frac{1}{m-1}}u(x,t) \to f(x) \quad\mbox{uniformly \,\,as}\,\,t\to+\infty.
          \end{equation}
          And  if $F$ is concave,  $z(x,t)$ converges to $f(x)$ in $ C_x^0(\overline\Omega)\cap C_x^{k,\al}(\Omega)$ for $k=1,2.$
\end{corollary}
\noindent {\sc Proof.}
 The first part (i) at Proposition \ref{prop-pme-app} directly gives the convergence of $z(x,t)$ to $f(x)$ as $t\to\infty$ uniformly in $\Omega.$  So we will see the second estimate.  For each compact subsets $K\Subset K'$ of $\Omega,$  uniform convergence implies 
$$\frac{1}{2}\inf_{K'}f \leq z(x,t)\leq 2\sup_{K'} f\quad \mbox{in}\,\,\,{K'\times[T,\infty)} $$
for large $T>1.$ For $w:=u^m,$ we have 
$$\frac{1}{2}\inf_{K'}f^m \, t^{-\frac{m}{m-1}} \leq w(x,t)\leq 2\sup_{K'} f^m\,t^{-\frac{m}{m-1}}\quad \mbox{in}\,\,\,{K'\times[T,\infty)}. $$
Let $t_o>2T.$ Then  there exist uniform constants $C_1, C_2$ with respect to time such that
$$C_1 t_o^{-\frac{m}{m-1}} \leq w(x,t)\leq C_2\,t_o^{-\frac{m}{m-1}}\quad \mbox{on}\,\,\,K'\times\left[\frac{t_o}{2},t_o\right].$$
We define $\tilde w(x,t)= t_o^{\frac{m}{m-1}}w(x,t_ot)$  in $K'\times\left[\frac{1}{2},1\right] $ and  we have $$C_1\leq \tilde w(x,t)\leq C_2 \quad \mbox{on}\,\,\,K'\times\left[\frac{1}{2},1\right].$$ Then $\tilde w$ 
 satisfies uniformly parabolic equation, $m\tilde w^{1-1/m}F(D^2\tilde w)=\tilde w_t$, in $K'\times(1/2,1]$  using scaling property.  
 Thus we get 
$$\|\tilde w(\cdot,1)\|_{C^{k,\al}_x(K)}\leq C\|\tilde w\|_{L^{\infty}(K'\times[1/2,1])}=C\| t_o^{\frac{m}{m-1}}w\|_{L^{\infty}(K'\times[t_o/2,t_o])}\leq C\|f^m\|_{L^{\infty}(\Omega)}$$
from the concavity  of $F,$
which means that  for any $t_o>2T,$
$$\| t_o^{\frac{m}{m-1}}w(\cdot,t_o)\|_{C^{k,\al}_x(K)}=\|\tilde w(\cdot,1)\|_{C^{k,\al}_x(K)} \leq C\|f^m\|_{L^{\infty}(\Omega)}.$$
Therefore uniform  convergence of $z^m$ to $f^m$ and uniform $C^{k,\al}_x$ estimates will give  that $z^m$ converges to $f^m$ in $C_x^{k,\al}$ - norm.
\qed

\subsection{ Square-root concavity of the pressure}
 Let  $v=u^{m-1}$ be the pressure and let  $v = w^2. $ 
We are going to prove the concavity of $w$ in spatial variables for $m>1.$   
The fact that $w$ is a suitable function to perform geometrical investigations
was demonstrated by Daskalopoulos, Hamilton and Lee  at \cite{DHL}.  We remark that the following computation is also valid for the fast diffusion, $\textsf{m}_{\Omega,F}<m<1$.

 First, let us approximate the equation: for $0<\eta<1,$
\begin{equation}
\begin{cases}
u_{\eta,t} = F(D^2u_{\eta}^m) \,\,\mbox{in}\,\, \Omega\times(0,\infty) \\
u_{\eta} = \eta \,\,\mbox{on}\,\, \partial\Omega\times(0,\infty)\\
u_{\eta,o} \geq \eta  \,\,\mbox{in}\,\, \Omega,
\end{cases}
\label{eq-pme-app-0}
\end{equation}
where we assume $\eta+\frac{1}{2}u_o<u_{\eta,o}\leq\eta+2 u_o.$
Let  $g_{\eta}=u_{\eta}^m$. Then $g_{\eta}$ satisfies the following equations:
\begin{equation}
\begin{cases}
mg_{\eta}^{1-1/m}F(D^2 g_{\eta})=g_{\eta,t}\,\,\mbox{in}\,\, \Omega \times(0,\infty)\\
g_{\eta} = \eta^m \,\,\mbox{on}\,\, \partial\Omega\times(0,\infty)\\
g_{\eta,o} > \eta^m \,\,\mbox{on}\,\, \Omega,
\end{cases}
\label{eq-pme-app-0-1}
\end{equation}
which is uniformly parabolic for a fixed $\eta>0$ since $g_\eta\geq \eta^m$ from the Comparison  principle.
We also assume that $g_{\eta,o}\in C^{\infty}(\Omega)$ and  $\eta^m+\frac{1}{2}g_o \leq g_{\eta,o}=g_{\eta}(\cdot,0)\leq 2g_o+\eta^m $ in $\Omega .$
Then we have the following uniform estimate with respect to $\eta$ so it suffices to show the concavity of $w_{\eta}.$

\begin{lemma}\label{lem-d} 
Let $F$ satisfy $(F1)$ and $F(0)=0$ and let $g_o\in \cC(\overline\Omega).$  
For each $t>s>0,$
 there are uniform constants $0<c_0(t),\, c_1,\, C_0(t,s)<\infty$ independent of $\eta>0$ such that
$$0<c_0(t)<|\D_x g_{\eta}|<c_1\quad \mbox{on}\,\,\, \p\Omega\times(0,t]$$ and 
$$|\D_x g_{\eta}| <C_0(t,s)\quad \mbox{on}\,\,\, \overline\Omega\times[s,t].$$
\end{lemma}
\noindent {\sc Proof.}
We establich a subsolution and a supersolution  of \eqref{eq-pme-app-0}.
Let $\vp^-$ be the positive eigen-function with respect to the eigenvalue $\mu^->0$ for the Pucci's operator  $\cM^-$ from Theorem \ref{thm1}, that is, $\vp^->0$ solves 
\begin{equation*}
\begin{cases}
 -\cM^-(D^2 \vp^-(x)) =\mu^- \vp^-(x)\quad\text{in $\Omega,$}\\
\vp(x)^-=0\quad\text{on $\partial\Omega.$}
\end{cases}
\tag{{\bf EV}}
\end{equation*}  
 We may assume that $g_{\eta,0}\geq \vp^-+\eta^m$  by multiplying a positive constant since $g_o\in \cC(\overline\Omega)$ and since $\cM^-$ is positively homogeneous of degree one.  
Define $K:=\mu^-m(1+||\vp^-||_{\infty})^{\gamma}>0$ for $\gamma:=1-\frac{1}{m}>0$ and  $h(x,t):=\eta^m+\vp^-e^{-Kt}.$ Then we have 
\begin{align*}
mh^{\gamma}F(D^2h)-h_t &\geq mh^{\gamma}\cM^-(D^2h)-h_t\\
&=mh^\gamma e^{-Kt}\left\{\cM^-(D^2\vp^-)+\frac{K\vp^-}{m(\eta+\vp^-e^{-Kt})^{\gamma}}\right\} \geq 0,\end{align*}
$h=\eta^m$ on $\p\Omega$  and $h(\cdot,0)\geq\vp^-+\eta^m.$ Thus the Comparison principle gives that  $g_{\eta}\geq h=\eta^m+\vp^-e^{-Kt},$  where $K$ depends on the initial data $g_o.$ So it follows that
$$|\D g_\eta(\cdot,t)|\geq c_o e^{-Kt}>0\quad \mbox{on}\,\,\, \p\Omega.$$ 

On the other hand,  let $\vp^+$ be the positive eigenfunction of 
\begin{equation*}
\begin{cases}
 -\cM^+(D^2 \vp^+(x)) = {\vp^{+}}^{\frac{1}{m}}(x)\quad\text{in $\Omega,$}\\
\vp(x)^+=0\quad\text{on $\partial\Omega.$}
\end{cases}
\tag{{\bf EV}}
\end{equation*} 
from Theorem \ref{thm1}. Multiplying a positive constant, we assume that $g_{\eta,0}\leq \vp^++\eta^m$ and $\vp^+$ is the  eigen-function with an eigen-value $\mu_o>0.$ If we define $h:=\vp^++\eta^m,$  then $h$ satisfies  $mh^\gamma F(D^2h)-h_t\leq mh^\gamma \cM^+(D^2h)-h_t\leq mh^\gamma(-\mu_o{\vp^+}^{1/m}) <0$ in $\Omega\times(0,\infty).$
From Comparison principle, we ontain  
$$g_{\eta}\leq \vp^++\eta^m,$$ which means that $$|\D g_{\eta}|< C_o=C_o(\vp^+)\quad\mbox{on}\,\,\,\p\Omega\times(0,\infty). $$ 
uniformly in $\eta>0.$   
A  similar argument as in (2),(ii) at Proposition \ref{prop-reg} gives 
$$|\D_xg_\eta|<C_o\quad\mbox{in}\quad\overline\Omega\times(s,t).$$
\qed

\begin{lemma}\label{lem-g-t}
 Suppose that $F$ satisfies $(F1),(F2), (F3).$  
Let $u$ be the solution of \eqref{eqn-pme} with $F(D^2u_o^m)\leq0$ in $\Omega$  and let $u_{\eta}$ be the solution of \eqref{eq-pme-app-0} with the initial data $u_{\eta,o} $ satisfying $F(D^2u^m_{\eta,o})\leq0.$ Then $u$ and $u_{\eta}$ are   nonincreasing in time.
\end{lemma}
\noindent {\sc Proof.}
According to Lemma \ref{lem-d}, it suffices to show that $g_{\eta,t}\leq0$ for any $\eta>0.$ Let us  fix $\eta>0$ and approximate the operator $F$ by smooth operators $\tilde F(\cdot):=F_{\e}(\cdot)-F_{\e}(0)$ in Lemma \ref{lem-F-approx}. Let $g_{\e,\eta}$ be the solution of \eqref{eq-pme-app-0-1} with the same initial data $g_{\eta,o}.$
For simplicity, we denote $g_{\e,\eta}$ and $\tilde F_{\e}$ by $g$ and $F,$ where the equation $\eqref{eq-pme-app-0-1}$ is uniformly parabolic in $\Omega\times(0,T]$ for a fixed  $\eta.$   Now define $$h:=g_t-\delta t-\delta$$ for small $\delta>0.$ Then $h$ is negative on the parabolic boundary. Indeed,  at $t=0$ we have 
  $h= mg^{1-1/m}F(D^2g)-\delta\leq mg^{1-1/m}\sqrt n\Lambda\e-\delta< 0$ for small $0<\e<<\delta$ and $h<0$ on $\p\Omega\times(0,T].$  Assume that there is $t_o\in(0,T]$ such that $h$ vanishes  at some point $x_o\in\Omega$ for the first time. Then   at the maximum point $(x_o,t_o),$ we have
\begin{align*}
0&\geq mg^{1-1/m}F_{ij}(D^2g)h_{ij}-h_t =-\left(1-\frac{1}{m}\right)\frac{g_t^2}{g} +\delta\\
&\geq-\left(1-\frac{1}{m}\right)\frac{\delta^2(t_o+1)^2}{\eta^m} +\delta\geq -\left(1-\frac{1}{m}\right)\frac{\delta^2(T+1)^2}{\eta^m} +\delta,
\end{align*}
which is  a  contradiction if we select $\delta$ and $\e$  small enough.
 Thus for a given $\eta, T>0,$  there is $\delta(\eta,T), \e(\eta,T)>0$ such that  if  $0<\delta<\delta(\eta,T)$ and  $0<\e<\e(\eta,T),$ then 
 $$\mbox{ $g_{\e,\eta,t}<\delta t+\delta\,\,$   and }\,\, g_{\e,\eta,t}\leq g_{\e,\eta}\,\,\,\,\,\mbox{ in $\,\,\Omega\times(0,T]$ }.$$
Letting $\e>0$ and $\delta>0$ go to $0,$ we have 
  $g_{\eta,t}\leq0$ in $\Omega\times(0,T]$  from the uniform Lipschitz estimates of $g_{\e,\eta}$ for a given $\eta>0.$ This completes the proof.
 \qed

\begin{lemma}\label{bdry} 
 Suppose that $F$ satisfies $(F1),(F2), (F3)$ 
and  that $\Omega$ is strictly convex. Let $u$ and $u_\eta$ be the solutions in Lemma \ref{lem-g-t}. Then  for each $T>0,$  there is $\eta(T)>0$ such that for $0<\eta<\eta(T),$ we have 
\begin{equation}\label{eq-bdry}
w_{\eta,\al\al}(x,t) = \frac{m-1}{2mg_{\eta}^{2-\frac{m-1}{2m}}}\left(g_{\eta} g_{\eta,\al\al} -\frac{m+1}{2m}g_{\eta,\al}^2\right) \leq
\,\mbox{sign}\,(1-m)\frac{c_o}{\eta^{\frac{m+1}{2}}}
\end{equation}
 on $(x,t) \in \partial\Omega\times(0,T]$  for any direction $e_{\al},$
where  $c_o>0\,$
is independent of $\eta>0.$
\end{lemma}
\noindent {\sc Proof.}
(a) Let us fix $\eta>0.$
 First, let us approximate the operator $F$ by {$F_{\e}(\cdot)-F_{\e}(0)$} as in Lemma \ref{lem-F-approx}. 
and  consider the approximated  equation: 
\begin{equation}
\begin{cases}
u_{\e,t} = F_{\e}(D^2u_{\e}^m) -F_{\e}(0)\,\,\mbox{in}\,\, \Omega\times(0,T], \\
u_{\e} = \eta \,\,\mbox{on}\,\, \partial\Omega\times(0,T],\\
u_{o,\e} > \eta \,\,\mbox{on}\,\, \Omega.
\end{cases}
\label{eq-pme-app}
\end{equation}
Let $g=u^m$ and $g_{\e}= u_{\e}^{m}.$
Then $g_{\e}$  satisfies
\begin{equation}
mg_{\e}^{\gamma}(F_{\e}(D^2g_{\e})-F_{\e}(0))= g_{\e,t} \quad\mbox{in}\,\, \,\Omega\times(0,T], \,\quad\gamma:=1-1/m>0,\\
\label{eq-pme-app-1}
\end{equation}
We will denote  $g_{\e},$ $F_{\e}$ by $g$, $F,$ respectively, for the simplicity. 

 Let us fix a boundary point $(x_o,t_o)\in\p\Omega\times(0, T].$ We denote $x_o$ by origin. {Now we introduce the coordinate system such that $x_0=0$
and that the tangent plane is $x_n=0$ at the origin.}
When $\tau=e_{i} ,(i=1,\cdots,n-1)$ is a tangential direction at $x_o = 0,$ $ g_{\tau}=0 $
and $g_{\tau\tau} = g_{\nu}\gamma_{\tau}$ at the origin where $e_\nu$ is  the outer normal  vector to $\p\Omega$
and $\gamma_{\tau}$ is the curvature of $\p\Omega$ in the direction $\tau.$ 

 (i.) According to boundary estimates at   Lemma \ref{lem-d}  and the  strict convexity of $\p\Omega,$ we have     $0<c(T)<-g_{\tau\tau} < C$ for any tangential vector $e_\tau$   and hence 
$$gg_{\tau\tau}(0,t)-\frac{m+1}{2m}g_{\tau}^2(0,t) \leq -c_o\eta^m\quad\mbox{on}\,\,\,\p\Omega\times(0,T]$$
for some $c_o(T)>0.$
 We also have $|g_{e_i,e_j}| \leq C$ on $\p\Omega,$ $ \,(1\leq i,j \leq n-1)$
 for some constant $C$ depending on $\p\Omega$ and  $C_o$ (which is a uniform bound for gradients on the boundary).

 (ii.)
        Near the origin, $\p\Omega$ is represented by $x_n=\gamma (x')=\frac12 B_{ij}x_ix_j+O(|x'|^3).$ The estimate 
        $c_0<\gamma_{\tau\tau}<C_0$  says that the eigen values of $\left(B_{ij}\right)$ is in $\left[c_0,C_0\right].$
        After a change of coordinate of {$\re^{n-1},$}
        the boundary becomes $x_n=\tilde{\gamma}(x')=\frac12 |x'|^2+O(|x'|^3)$ and  the operator $F$ will be transformed to a new operator $\tilde{F}$
        with new elliptic coefficients $\tilde{\lambda}=\tilde{\lambda}(\lambda,\Lambda,c_0,C_0)$ and        
        $\tilde{\Lambda}=\tilde{\Lambda}(\lambda,\Lambda,c_0,C_0)$ that are uniformly bounded and positive.
So $\p\Omega$ is close to  a unit ball with an error $O(|x'|^3)$ near the origin. 
 For simplicity we are going to  assume that $\Omega=B_1(e_n).$ The general domain can be considered with a simple modification as \cite{CNS}.

(iii.)  We claim that $|g_{e_i,e_n}(0,t_o)| \leq C$ for $1\leq i \leq n-1.$ \\
 For positive constants $A,B$ and $D$, let us define in $\Omega\times(t_o/2,t_o)$
\begin{align*}
w_{\pm}(x,t) &= \left\{\p_{T_k}g \pm A\sum_{l=1}^{n-1}g_l^2\pm D x_n^2\right\}T(t)\,
= \left\{(1-x_n)g_k+x_kg_n \pm A\sum_{l=1}^{n-1}g_l^2\pm Dx_n^2,
\right\}T(t)\end{align*}
where $\p_{T_k}g := (1-x_n)g_k+x_kg_n  $ is a directional derivative and coincides with
a tangential derivative on $\p B_1$ and $T(t):=e^{M(t-t_o/2)}-1.$ Let $v:=C(A+D)x_nS(t)$ for $S(t):=K(t-t_o/2)  \geq 0.$ The constant $M,K>0$   will be chosen so that  $S(t)\geq T(t) $ in $(t_o/2,t_o).$

Since $g=\eta^m$ on $\p B_1$ and 
\begin{align*}
\frac{1}{2}g_l^2 &= \frac{1}{2}[(1-x_n)g_l+x_lg_n+x_ng_l-x_lg_n ]^2 \leq [(1-x_n)g_l+x_lg_n]^2+[x_ng_l-x_lg_n ]^2\\
&\leq [(1-x_n)g_l+x_lg_n]^2+C|x|^2 \,\,\,\,\,\,\,(\mbox{we recall}\,\,|\nabla g|<C \,\,\mbox{on}\,\,\overline  B_1\times(t_o/2,t_o)),
\end{align*}
we see that for all $x\in\p{B_1}$
\begin{align*}
-C(A+D)|x|^2 T(t)\leq w_-\leq w_+ \leq C(A+D)|x|^2 T(t).
\end{align*}
Since $|x|^2=2x_n$ for any $x\in\p B_1,$ we obtain that
$$-C(A+D)x_nT(t)\leq w_-\leq w_+ \leq C(A+D)x_nT(t)
\quad\mbox{for }\,\,\,x\in\p_p( B_1\times(t_o/2,t_o)) $$ 
and then we have
$$-v\leq w_-\leq w_+ \leq v
\quad\mbox{for any}\,\,\,x\in\p_p( B_1\times(t_o/2,t_o)). $$ 
Now, let us consider a linearized operator 
$$H[u] = mg^{\gamma}F_{ij}(D^2g)D_{ij}u-u_t.$$
If $H[w_+]\geq H[v]$ and $H[w_-]\leq H[-v]$ in $B_1\times(t_o/2,t_o) $ for some constants $A$  $D,$ $M$ and $K,$   then the   comparison principle gives 
$$ -v \leq w_-\leq w_+ \leq v \,\,\,\mbox{in}\,\,\,B_1\times(t_o/2,t_o).$$
Therefore, we deduce that, for $1\leq k\leq n-1,$
\begin{align*}
|g_{kn}(0,t_o)|&=|(w_+)_n(0,t_o)|\leq C(A+D)S(t_o).
\end{align*}
So, it remains to  show that $H[w_+] \geq H(v)$ and $H[w_-]\leq H(-v)$ if $A$(uniform with respect to $\eta,\e$) and $D$ are chosen    large enough.
Using $mg^{\gamma}F(D^2g)- g_t=0$ and the ellipticity of $F,$ it follows that
\begin{align}\label{eq-f}
&|mg^{\gamma}g_{nn}|^2 \leq C\left(|mg^{\gamma}|^2\sum_{(i,j)\neq(n,n)}|g_{ij}|^2 +|g_t|^2\right)\,\,\,\mbox{in}\,\,\Omega,\\
& mg^{\gamma}\sum_{(i,j)\neq(n,n)}|g_{ij}|^2 \geq  Cmg^{\gamma}g_{nn}^2- C\frac{|g_t|^2}{g^{\gamma}}\,\,\mbox{in}\,\,\Omega.
\end{align}
Using the above inequalities, we have 
\begin{align*}
H[w_+]&=-T(t)\gamma\frac{g_t}{g}\left\{(1-x_n)g_k+x_kg_n+ 2A\sum_{l=1}^{n-1}g_l^2\right\}\\
&+T(t)2mg^{\gamma}\left\{-\sum_{i=1}^{n}F_{ni}g_{ki}+\sum_{i=1}^{n}F_{ki}g_{ni} +A\sum_{l=1}^{n-1}F_{ij}g_{li}g_{lj} +DF_{nn}\right\}\\
&-T'(t)\left\{(1-x_n)g_k+x_kg_n + A\sum_{l=1}^{n-1}g_l^2 +Dx_n^2
\right\}\\
&\geq -T(t)|\gamma|\frac{|g_t|}{g}\left|(1-x_n)g_k+x_kg_n+ 2A\sum_{l=1}^{n-1}g_l^2\right|+T(t)D\lambda m\eta^{m\gamma}\\
&+2T(t)m g^{\gamma}\left\{ -C(\sum_{l,i=1}^{n}|g_{li}|^2)^{\frac{1}{2}} +\frac{\lambda A}{C}\sum_{l,i=1}^{n}g_{li}^2+A\lambda/2\right\}-T(t)C\displaystyle\frac{|g_t|^2}{|g|^2}g^{2-
\gamma}\\
&-T'(t)\left\{(1-x_n)g_k+x_kg_n + A\sum_{l=1}^{n-1}g_l^2 +Dx_n^2
\right\}\\
&\geq 2T(t)m g^{\gamma}\left\{ -C(\sum_{l,i=1}^{n}|g_{li}|^2)^{\frac{1}{2}} +\frac{\lambda A}{C}\sum_{l,i=1}^{n}g_{li}^2+A\lambda/2\right\}\\
&+CT(t)\left(D\eta^{m\gamma}-C\right)-CT'(t)\left(D-C
\right).
\end{align*}
We note that $g, \displaystyle\left|\frac{g_t}{g}\right|$ and $|\D g|$ are  uniformly bounded with respect  to $\eta $ and small $\e(\eta)$ in $B_1\times(t_o/2,t_o)$ according to Aronson-Benilan inequality at Lemma \ref{lem-ab} and Lemma \ref{lem-d}, \ref{lem-g-t}. 
Thus if $A>C/\sqrt{\lambda},\,\, D>2A$ and { if $D\min(\eta^{m\gamma},1)$} is big enough, we get 
\begin{equation*}
H[w_+]\geq   C_1D(C_2\eta^{m\gamma}T-T'). 
\end{equation*}
Setting $M=C_2\eta^{m\gamma}=C_2\eta^{m-1}$ and $K=\frac{2}{t_o}\left(e^{C_2\eta^{m-1}t_o/2}-1\right),$ 
we have 
$$H[w_+]\geq 0 \geq H[v]=-C(A+D)x_nS'(t).$$
Similarly, we have 
$H[w_-]\leq H[-v]$ in $B_1\times(t_o/2,t_o). $ 
Therefore, we have proved that
\begin{align*}
|g_{kn}(0,t_o)|&=|(w_+)_n(0,t_o)|\leq \frac{1}{\eta^{m-1}}\left(e^{C_2\eta^{m-1}t_o/2}-1\right)\leq  2C_2 T,
\end{align*}
for $1\leq k\leq n-1$ and 
for small $0<\eta<\eta(T),$ where $C_2$ and $\eta(T)$ are uniform with respect to $\eta,\e.$

(iv) Lastly, since $g^2_{nn}(0,t_o)\leq C\sum_{(i,j)\neq(n,n)}|g_{ij}|^2 $ from $\eqref{eq-f},$  we have   $$|g_{nn}|(0,t_o)\leq  C (T )\,\,\,\,
\mbox{and}\,\,\,\,|D^2 g|(0,t_o)\leq  C( T) ,$$
where $C(T)$ is independent  of $\eta>0$ and $\e>0.$ 
Therefore, we have that for any unit vector $e_{\beta}:=\beta_1e_{\tau}+\beta_2e_{\nu},$
\begin{align*}
gg_{\beta\beta}(0,t_o)-\frac{m+1}{2m}g_{\beta}^2(0,t_o)& \leq   g\left(-c_o\beta_1^2+C(T)\left(\beta_2^2+2\beta_1\beta_2\right)\right)-\beta_2^2\delta_o\\
&\leq-\frac{c_o}{2}\eta^m\beta_1^2+\left(\eta^mC(T)\left\{\frac{1}{2c_o}+1\right\}-\delta_o\right)\beta_2^2\,\,\leq -\frac{c_o}{2}\eta^m
\end{align*}
for a small $\eta>0, $ using Young's inequality and gradient estimate at Lemma \ref{lem-d}. 

(b) For the general operator instead of smooth operators, the result follows from the uniform $C^{2,\be}$- estimates  since $g_{\e}$ satisfies the uniformly parabolic equation with ellipticity constants  related to $\eta>0,$ \eqref{eq-pme-app-1}. 
\qed

\begin{remark}{\rm
\item
\begin{enumerate}[(i)]
\item The boundary estimate \eqref{eq-bdry} holds if $|D^2g_{\eta}|$ is uniformly bounded  in $\overline\Omega\times(0,T]$ with respect to $\eta>0.$ In Lemma \ref{bdry},
 we have proved the estimate for the solutions with initial condition that $F(D^2g_{\eta,o})\leq 0.$
  \item To prove the estimate $\eqref{eq-bdry}$ up to the boundary, we need to prove a weighted $C^{2,\al}_{\delta}$- estimate of $u_{\eta}=g_{\eta}^{1/m}$ up to the boundary, 
which will be studied in the future work. When $F(D^2u)=\La u,$
Schauder theory has been proved in \cite{KL}.
\end{enumerate}
}
\end{remark}

\begin{lemma}\label{lem-c2f-p}
{ Let $F\,$ satisfy (F1), (F2) and (F3)} and let $\Omega\,$ be a strictly  convex
bounded domain. Let $u$ be the solution of $\eqref{eqn-pme}$ with an initial data $u_o^m\in\cC(\overline\Omega).$
and let $u_{\eta}$ be an approximated solution of \eqref{eq-pme-app-0-1}. 
Assume  that the boundary estimate of Lemma \ref{bdry}  holds for approximated solutions $u_{\eta}.$

 If
$\sqrt{v_o}=u_o^{\frac{m-1}{2}}\,$ is concave, then the pressure  $v(x,t)=u^{m-1}(x,t)$ of
$\eqref{eq-pres}$ is square root-concave in the spatial
variables, i.e., $D_x^2\sqrt{v(x,t)} \leq 0$ in $\Omega\times(0,\infty).$
\end{lemma}

\noindent {\sc Proof.}
(i) First, we fix $T>0.${ We may assume that $u_{\eta,o}\in C^{\infty}(\Omega)$ satisfies $\eta^m+\frac{1}{2}u_o^m\leq u^m_{\eta,o}\leq \eta^m+2u^m_o, \,\,F(D^2u_{\eta,o}^m)\leq0$ and  $D^2\sqrt{u^{m-1}_{\eta,o}} \le 0$  in $\Omega $  } 
and also assume that 
there is small $\eta(T)$ such that for $0<\eta<\eta(T),$ the boundary estimate \eqref{eq-bdry} is true from the assumption.  

If we show the square  root - concavity of $v_{\eta}$, the pressure  of $u_{\eta}$   in $\Omega\times(0,T],$ then the concavity for $v=u^{m-1}$ follows from uniform convergence.   Indeed, the uniform Lipschitz estimates of $u_{\eta}^{m}(x,t)$ will give us  uniform convergence of $u_{\eta}\,\,$
to $u$ in each compact subset of $\Omega\times(0,T],$ from Lemma \ref{lem-d}, and the limit $u$  also satisfies
$$u^{\frac{m-1}{2}}(x,t)+u^{\frac{m-1}{2}}(y,t)-2u^{\frac{m-1}{2}}\left(\frac{x+y}{2},t\right) \leq 0.$$

(ii) Now, 
let us fix $T$ and  $\eta$ for $0<\eta<\eta(T).$   It remains to show that $u_{\eta}^{\frac{m-1}{2}}$ is concave in $\Omega\times(0,T]$ for small $\eta>0.$ 
Let us   {approximate $F$ by a smooth $F_{\e}$ from  Lemma \ref{lem-F-approx}} 
and let $u_{\e,\eta}$ be the solution of the approximated equation
\begin{equation}
\begin{cases}
u_{t} = F_{\e}(D^2u^m) -F_{\e}(0)\,\,&\mbox{in}\,\, \Omega\times(0,T], \\
u = \eta \,\,&\mbox{on}\,\, \partial\Omega\times(0,T],\\
u(\cdot,0)=u_{\eta,o} > \eta \,\,&\mbox{on}\,\, \Omega.
\end{cases}
\label{eq-pme-app}
\end{equation}
For simplicity, we denote $u_{\e,\eta}, u_{\e,\eta}^m$ by $u, g=u^m.$ The function $g$ solves $$ mg^{\gamma}(F_{\e}(D^2g)-F_{\e}(0))= g_t \quad\mbox{in}\,\, \Omega\times(0,T], \,\quad(\gamma=1-1/m>0),$$
which is uniformly parabolic for a given $\eta>0.$ 
 
The geometric quantity $w:=\sqrt{v}=u^{\frac{m-1}{2}}$ satisfies
$$ w_t =\frac{m-1}{2}w^{\frac{m-3}{m-1}} F_{\e}\left(\frac{2m}{m-1}w^{\frac{3-m}{m-1}}\left(w^2D^2w + \frac{m+1}{m-1}wDwDw^T\right)\right)-\frac{m-1}{2}w^{\frac{m-3}{m-1}} F_{\e}(0).$$ 
After the change of the time $t \mapsto mt,\,$ the equation will be simplified into
\begin{equation}
w_t = \frac{m-1}{2m}w^{\frac{m-3}{m-1}} F\left(\frac{2m}{m-1}w^{\frac{3-m}{m-1}}\left(w^2D^2w +  rwDwDw^T\right)\right)-{\frac{m-1}{2m}w^{\frac{m-3}{m-1}} F_{\e}(0)}\end{equation}
with $ r =\frac{m+1}{m-1}.$ 
By taking differentiation, we have
\begin{align*}
w_{\al\be t} &= \frac{m-1}{2m}w^{\frac{m-3}{m-1}}F_{ij,kl}\cdot\left(\frac{2m}{m-1}w^{\frac{3-m}{m-1}}\left(w^2D_{ij}w + \frac{m+1}{m-1}wD_iwD_jw\right)\right)_{\alpha}\\
&\cdot\left(\frac{2m}{m-1}w^{\frac{3-m}{m-1}}\left(w^2D_{kl}w + \frac{m+1}{m-1}wD_{k}wD_lw\right)\right)_{\be}\\
&+ F_{ij}\cdot(2w_{\al}w_{\be}D_{ij}w+2ww_{\al\be}D_{ij}w+2ww_{\al}D_{ij}w_{\be}+2ww_{\be}D_{ij}w_{\al}+w^2D_{ij}w_{\al\be}\\
&+rw_{\al\be}D_iwD_jw+2rw_{\al}D_iw_{\be}D_jw+2rw_{\be}D_iw_{\al}D_jw\\
&+2rwD_iw_{\al}D_jw_{\be}+2rwD_iw_{\al\be}D_jw)\\
& +\frac{m-1}{2m}\frac{m-3}{m-1}w^{\frac{m-3}{m-1}-1}w_{\al \be} F\left(\frac{2m}{m-1}w^{\frac{3-m}{m-1}}\left(w^2D^2w +  rwDwDw^T\right)\right)\\
&-\frac{m-3}{m-1}w^{-1}w_{\al \be} F_{ij}\left(w^2D_{ij}w +  rwD_iwD_jw\right)\\
&-\frac{m-1}{2m}\frac{m-3}{m-1}\frac{2}{m-1}w^{\frac{m-3}{m-1}-2}w_{\al}w_{\be} F\left(\frac{2m}{m-1}w^{\frac{3-m}{m-1}}\left(w^2D^2w +  rwDwDw^T\right)\right)\\
&+\frac{m-3}{m-1}\frac{2}{m-1}w^{-2}w_{\al}w_{\be} F_{ij}\left(w^2D_{ij}w +  rwD_iwD_jw\right)\\
&-\frac{m-1}{2m}\frac{m-3}{m-1}w^{\frac{m-3}{m-1}-2}\left\{ww_{\al \be}-\frac{2}{m-1}w_{\al}w_{\be} \right\} F_{\e}(0),
\end{align*}
for $F_{ij}= F_{\e,ij}\left(\frac{2m}{m-1}w^{\frac{3-m}{m-2}}\left(w^2D^2w +  rwDwDw^T\right)\right).$ 

In order to show the concavity of $w,$ consider 
$$
\sup_{y\in \Omega }\sup_{|e_{\beta}|=1} \, w_{\beta\beta}(y,t)+\psi(t),
$$
 where  $e_{\beta}\in S^{n-1}$ and 
a negative  function $\psi(t)$  with $\psi(0)<0$ will be chosen later.
Let us assume that
$$
\sup_{y\in \Omega }\sup_{|e_{\beta}|=1} \, w_{\beta\beta}(y,t)+\psi(t)=0\quad\mbox{at}\,\,\,t=t_o,
$$
for the first time.
From the assumption  that the pressure is  initially  square-root concave, the quantity   $\displaystyle\sup_{y\in \Omega }\sup_{|e_{\beta}|=1} \, g_{\beta\beta}(y,t)+\psi(t)$ is negative at $t=0.$ 

Now, we assume that the supremum
$$\sup_{y\in \Omega}\sup_{|e_{\beta}|=1}  w_{\be\be}(x,t_o) = w_{\oa\oa}(x_o,t_o) = -\psi(t_o)(>0)$$
is achieved at $(x_o,t_o)\in \overline\Omega\times(0,T]$ with a unit vector  $e_{\oa} $ and    assume that $x_o = 0$ without losing of generality.   Then, 
the assumption on the boundary that $w_{\eta,\be\be}\leq 0$
yields that $(0,t_o)\,$ should be an interior point. 
We introduce an orthonormal coordinates in which $e_{\oa}$ is taken as one of the coordinate axes and we assume that 
$$w_{\oa\be}(0,t_o)=0\quad\mbox{if}\,\,\,\be\neq\oa.$$

In order to create extra terms,  we perturb second derivatives of $w$ and we   
 use the function
$$Z(x,t) = w_{\al\be}(x,t)\xi^{\al}(x)\xi^{\be}(x)$$
where  $\xi^{\be}(x)=\delta_{\oa\be} + c_{\oa}x^{\be} + \frac{1}{2}c_{\oa}c_{\gamma}x^{\gamma}x^{\be} .$  We are going to choose $c_{\al}\,$ so that
$-4w^2c_{\al}+4ww_{\al}=0\,\,$ at the maximum point $(x,t)=(0,t_o)$ an then the function $Z$ will help the third derivatives   cancel out, which   appear  in the  porous medium equation after differentiations.   
We note  that at the maximum point $(0,t_o),$ we have 
$$w_{\oa,\be} = 0\,\,\mbox{if}\,\,\be\neq\oa, \quad D_x^2Z \leq 0,\,\,\,\mbox{and}\,\,\,\D_x Z = 0$$
since  
\begin{align*}
D^2w(x,t)&<-\psi(t){\bold{I}}\quad\mbox{for}\,\,\,0<t<t_o,\\
Z(x,t)&=\displaystyle\vec \xi^{\,\, T} D^2w\, \vec \xi,\qquad 
Z(0,t_o)=w_{\oa\oa}(0,t_o)=-\psi(t_o).
\end{align*}
A simple computation gives us, at $(x,t)=(0,t_o),$
\begin{align*}
Z_i &= w_{\al\be i}\xi^{\al}\xi^{\be} + 2w_{\be i}c_{\al}\xi^{\al}\xi^{\be}\\
Z_{ij} &= w_{\al\be ij}\xi^{\al}\xi^{\be} + 4w_{\be ij}c_{\al}\xi^{\al}\xi^{\be}
+ 2w_{\be i}c_jc_{\al}\xi^{\al}\xi^{\be}
+ 2w_{ij}c_{\al}c_{\be}\xi^{\al}\xi^{\be}.
\end{align*}
and hence we have at $(0,t_o),$
\begin{align*}
Z_t &= w_{\al\be t}\eta^{\al}\eta^{\be}\\
&\leq w^2F_{ij}\cdot Z_{ij} +F_{ij}w_{\be
ij}\eta^{\al}\eta^{\be}(-4w^2c_{\al}+4ww_{\al}) +2F_{ij}w_{ij}(-w^2c_{\oa}^2+w_{\oa}^2) \\
&+(2wF_{ij}D_{ij}w+rF_{ij}D_iwD_jw)w_{\oa\oa} +4rF_{\oa
j}D_jww_{\oa}w_{\oa\oa} -2w^2F_{\oa
j}c_jc_{\oa}w_{\oa\oa }  \\
&+2rwF_{ij}D_jwD_iw_{\oa\oa} +rwF_{\oa\oa}w_{\oa\oa}^2\\
&+\frac{|m-3|}{2m}w^{-\frac{m+1}{m-1}}\left\{ww_{\oa\oa}+\frac{2}{m-1}w_{\oa}^2 \right\} C\e,\\
&\leq (2wF_{ij}D_{ij}w+rF_{ij}D_iwD_jw)Z+4rF_{\oa
j}D_jww_{\oa}Z -2w^2F_{\oa
j}c_jc_{\oa}Z  \\
&+rwF_{\oa\oa}Z^2 +\frac{|m-3|}{2m}w^{-\frac{m+1}{m-1}}\left\{ww_{\oa\oa}+\frac{2}{m-1}w_{\oa}^2 \right\} C\e
\\
&\leq (2wn\Lambda(-\psi)+Cr\Lambda|\D w|^2)Z+rwF_{\oa\oa}Z^2 +\frac{|m-3|}{2m}w^{-\frac{m+1}{m-1}}\left\{-w\psi+\frac{2}{m-1}w_{\oa}^2 \right\} C\e.
\end{align*}

Now let us define {$$Y(x,t):= Z(x,t)+\psi(t).$$}
We notice  that $Y(x,0)<0$ for any $x\in\Omega$ since we know $$\displaystyle\sup_{x\in \Omega}\sup_{|e_{\beta}|=1}  w_{\be\be}(x,0) <0\quad\mbox{and}\quad\psi(0)<0$$ and    $\p_t Y(0,t_o)\geq0$ since  we have 
\begin{align*}
D^2w(0,t)&<-\psi(t){\bold{I}}\quad\mbox{for}\,\,\,0<t<t_o,\quad Z(0,t_o)=w_{\oa\oa}(0,t_o)=-\psi(t_o).
\end{align*}
Thus we obtain at the maximum point $(0,t_o),$
$$0\leq \p_tY(0,t)=Z_t+\psi_t\leq \psi_t+K(\psi^2-\psi+\e) ,$$
where $K=c(n,m,\Lambda)\displaystyle\sup_{\Omega\times(0,T]}\left((w+w^{-\frac{2}{m-1}})+(1+w^{-\frac{m+1}{m-1}})|\D w|^2+ w\right).$ We note that $K$ is independent of $\e$ (and $\delta$) since $\D_x g_{\e,\eta}$ is uniformly bounded  and   $\eta>0$ is given.

If we set $\psi(t):=-\e\displaystyle-e^{-1/\delta}e^{Kt}\tan(K\sqrt\delta t),$ a simple calculation yields  
$$
 \psi_t+K (\psi^2-\psi+\e) <0$$ for $0<\e<<\delta$ and $K \sqrt\delta t<\frac{\pi}{2},$ which
implies  a contradiction if $t_o<\frac{\pi}{2K\sqrt\delta}.$
 Therefore, we deduce that  $$\sup_{y\in \Omega }\sup_{|e_{\beta}|=1} \, w_{\beta\beta}(y,t) < -\psi(t)\,\,\,\,\mbox{ for}\,\,\,  0<t<\frac{\pi}{2K\sqrt\delta}\,\,\,\mbox{ and for small} \,\,\,\e<<\delta,$$ and hence
   $$\p_{\be\be} w=\p_{\be\be} u_{\e,\eta}^{\frac{m-1}{2}} <\,\e+  e^{-1/\delta +K /\sqrt\delta}\quad \mbox{for any unit} \,\,\,e_{\be}$$
for $\displaystyle0<t<\frac{\pi}{8 K\sqrt{\delta}}$ and $0<\e<<\delta.$
Letting $\delta\ra 0,$ we conclude that
 $$\p_{\be\be} u_{\eta}^{\frac{m-1}{2}} \leq 0 \quad\,\,\,\mbox{in}\,\,\,\Omega\times(0,T] \,\,\,\,\mbox{for any }\,\,\,e_{\be},$$ that implies 
$$u_{\eta}^{\frac{m-1}{2}}(x,t)+u_{\eta}^{\frac{m-1}{2}}(y,t)-2u_{\eta}^{\frac{m-1}{2}}\left(\frac{x+y}{2},t\right) \leq 0.$$
\qed

\begin{corollary}
Let us assume that $F$ satisfies $(F1)$,$(F2)$ and $(F3)$ and $\Omega$ is  convex. 
If  $\sqrt {v_o}$ is concave, so is the viscosity solution $v(x,t).$
\end{corollary}

\begin{corollary}[\bf Square-root Concavity]\label{lem-p}
Let   $F$ satisfiy (F1), (F2) and (F3). 
If $\Omega$ is convex,  then $f^{m-1}$ is
square root-concave,  where $f$ is the positive eigenfunction of \eqref{eq-main-1-pme}.
\end{corollary}
\noindent {\sc Proof.} First, we may assume that $F(\cdot)\leq trace(\cdot)$ after a simple transformation from (F1) and (F3). 
Let $\phi_{\La}$ be the eigenfunction of \eqref{eq-main-1-pme} for Laplacian and let $u$  the solution of $\eqref{eqn-pme}$ with an initial data $u_o:=\phi_{\La}.$  
Then $\phi_{\La}$ satisfies $F(D^2\phi^{m}_{\La})\leq \La\phi^m_{\La}=-\frac{1}{m-1}\phi_{\La}<0$ in $\Omega$ and  has nontrivial bounded gradient on $\p\Omega.$ Moreover,   $\phi^{\frac{m-1}{2}}_{\La}$ is concave in a convex domain from  \cite{Ka,LV2}. 
Thus Lemmas  \ref{bdry} and \ref{lem-c2f-p} imply  $$D_x^2\sqrt{u^{m-1}}\leq0\quad\mbox{in $\Omega\times(0,\infty).$}$$
The uniform convergence at Proposition \ref{prop-pme-app},  Corollary $\ref{cor-pme-app}$ that is,
 $$tu^{m-1}(x,t)\to f(x) \quad \mbox{uniformly \,in}\,\,\Omega \,\,\,\mbox{ as $t\to+\infty,$}$$   will  preserve  the concavity of $f^{\frac{m-1}{2}}.$   Therefore, it follows that
$$f^{\frac{m-1}{2}}(x)+f^{\frac{m-1}{2}}(y)-\frac{1}{2}f^{\frac{m-1}{2}}\left(\frac{x+y}{2}\right)\leq0\quad \mbox{for}\,\,\,x,y\in\Omega.$$
\qed 

Now we   state the  strict  concavity of  solutions to \eqref{eq-main-1-pme},\eqref{eqn-pme} which follow from \cite{LV2}.
\begin{lemma}[ Strict Square-root Concavity]\label{lem-sc-pme}
Suppose that  $F$ satisfies (F1), (F2) and (F3)  and is  differentiable.  If $\Omega\,$ is smooth and strictly
convex, $f^{\frac{m-1}{2}}(x)\,$ is strictly  concave: there exists a
constant $c_1 > 0\,$ such that
\begin{equation*}
D^2\sqrt{h(x)} \leq -c_1\bf{I}.
\end{equation*}
\end{lemma}

\begin{theorem}[Eventual square root-concavity]
\label{theor2}
 We assume the same hypothesis as Lemma \ref{lem-sc-pme}.  Let $u_o\in\cC(\overline\Omega).$ Then, the pressure
$v(x,t)=u^{m-1}(x,t)$   is strictly square root -concave in $x$  variables for  large $t>0$. More precisely, for any 
$\e>0,$ there is $t_0=t_0(u_o,\e)$ such that
\begin{equation*} 
 D^2\sqrt{tv(x,t)}\leq
-(c_1-\e)\,{\bf I}
\end{equation*}
for $t\geq t_0$ and $x\in\Omega_{\e}=\{x\in\Omega |
\,d(x,\p\Omega) > \e\},$ where $c_1>0$ is the constant in
Lemma {\rm \ref{lem-sc-pme}}.
\end{theorem}

\begin{remark}{\rm
\item
\begin{enumerate}[(i)]
\item 
$F(\cdot)$ in Theorem \ref{theor2} is basically  Laplacian after a simple transformation if $F(\cdot)$ is differentiable and satisfies  (F1), (F2).
\item Condition $(F2)$ is required to have the convergence of $t v(x,t)$ to $f^{m-1}(x)$ as $t\ra+\infty$ and 
 the concavity of $F$ is required when we consider a concavity of solutions.
\end{enumerate}
}
\end{remark}


{\bf Acknowledgement} Ki-Ahm Lee was supported by Basic Science Research Program through the National Research Foundation of Korea(NRF)
grant funded by the Korea government(MEST)(2010-0001985) .

\end{document}